\def\R{{\bf R}}

\def\C{{\bf C}}

\def\Z{\bf Z}
\def\N{{\bf N}}
\def\H{{\cal H}}

\def\E{{\cal E}}
\def\K{{\cal K}}
\def\L{{\cal L}}

\def\D{{\cal D}}
\def\A{{\cal A}}

\def\Arg{\mathop{\rm Arg}\nolimits}

\def\Im{\mathop{\rm Im}\nolimits}

\def\diag{\mathop{\rm diag}\nolimits}
\def\const{\mathop{\rm const}\nolimits}
\def\sbp{\subparagraph}

\documentstyle[12pt]{article}
\newtheorem{theorem}{Theorem}
\newtheorem{corollary}{Corollary}

\newtheorem{statement}{Statement}
\newtheorem{lemma}{Lemma}
\newtheorem{remark}{Remark}
\newtheorem{example}{Example}

\newtheorem{definition}{Definition}

\title{Geometrical Quantization in Fock Space}
\author{
   V.P.Maslov and O.Yu. Shvedov  \\
{\small{\em Sub-faculty of Quantum Statistics and Field Theory,  }}\\
 {\small{\em Department of Physics, Moscow State University, }}\\
{\small{\em Vorobievy gory, Moscow 119899, Russia}} }

\begin{document}

\maketitle

{\small
We investigate an infinite dimensional analog of the theory
of Lagrangian manifolds with complex germs. To such a manifold
we assign a canonical operator that depends on creation and
annihilation operators. This operator is by definition the
geometrical quantization for these isotropic manifolds with
complex germs. We prove that for
secondary quantized equations this quantization is
the asymptotics for the Cauchy problem. Results of Berezin
are used thouroughly in the construction of the canonical
operator and in proofs of the theorems.
}
\newpage

\section{Introduction}

Construction of asymptotic solutions for multiparticle
Schr\"odinger and Liouville equations as the number
of particles tends to infinity was investigated in  \cite{MS1,MS2}.
The construction looks like follows. Both Schr\" odinger
and Liouville multiparticle equations may be presented in
a unified form through the creation and annihilation operators:

\begin{equation}\label{e1}
i\frac{\partial \hat\Phi}{\partial t}
=\left( T_{mn}\hat\psi_m^+\hat\psi_n^- +
\frac{\epsilon}{2}V_{klrs}\hat\psi_k^+\hat\psi_l^+
\hat\psi_r^-\hat\psi_s^-\right)\hat\Phi
\end{equation}

Here  $\hat \Phi$ is a Fock space element, $\hat\psi_j^{\pm}$
are creation and annihilation operators in this space \cite{B,LL},
 and we sum over
the repeating indices, $m,n,k,l,r,s=\overline{1,\infty}$.

The coefficients $T_{mn}$ and $V_{klrs}$ are

\begin{enumerate}

\item for the Schr\" odinger equation

\begin{eqnarray*}
T_{mn} & = & \int dx f_m^*(x)(-\Delta/2+U(x))f_n(x), \\
V_{klrs} & = & \int dxdyf_k^*(x)f_l^*(y)V(x,y)f_r(y)f_s(x), \\
& & x,y\in\R^{\nu}, \nu\in\N,
\end{eqnarray*}
where $\{f_1,f_2,\dots\}$ is an orthonormal basis in
$L^2(\R^{\nu})$, $\Delta$ is the Laplace operator in $\R^{\nu}$,
$U$ is an external potential, $V$ is the potential of the
interparticle interaction. Form (\ref{e1}) of the multiparticle
Schr\" odinger equation originates to papers  \cite{D,F1,F2,JW}.

\item for the Liouville equation

\begin{eqnarray*}
T_{mn} & = & i \int dpdqf_m^*(p,q)\left(\frac{\partial U}{\partial q}(q)
\frac{\partial}{\partial p} - p\frac{\partial}{\partial q}\right)
f_n(p,q),\\
V_{klrs} & = & i \int dp_1dp_2dq_1dq_2 f_k^*(p_1,q_1)f_l^*(p_2,q_2)\\
& & \cdot \left(\frac{\partial V(q_1,q_2)}{\partial q_1}
\frac{\partial}{\partial p_1} + \frac{\partial V(q_1,q_2)}
{\partial q_2}\frac{\partial}{\partial p_2}\right)
f_r(p_1,q_1)f_s(p_2,q_2), \\
& & p,q,p_1,q_1,p_2,q_2\in\R^{\nu}, \nu\in\N,
\end{eqnarray*}
where $\{f_1,f_2,\dots\}$ is an orthonormal basis in $L^2(\R^
{2\nu})$, $U$ and $V$ are as in the previous case an
external potential and the interparticle interaction potential
respectively. Sch\" onberg \cite{Sn1,Sn2} was the first one who suggested
the form (\ref{e1}) for the multiparticle Liouville equation,
see also  \cite{MT}.
\end{enumerate}

The asymptotics of the solution of equation (\ref{e1}) was
constructed as follows. After the substitution
$\sqrt{\epsilon}\hat\psi_j^{\pm}=\hat\phi_j^{\pm}$
equation (\ref{e1}) becomes

\begin{equation}\label{e2}
i\epsilon\frac{\partial\hat\Phi}{\partial t} = H
(\hat\phi_j^+,\hat\phi_j^-)\hat\Phi,
\end{equation}
for

$$ H(\hat\phi_j^+,\hat\phi_j^-) = T_{mn}\hat\phi_m^+\hat\phi_n^-
+\frac{1}{2}V_{klrs}\hat\phi_k^+\hat\phi_l^+\hat\phi_r^-\hat\phi_s^-$$
where
\begin{equation}\label{e3}
[\hat\phi_j^-,\hat\phi_l^+] = \epsilon\delta_{jl}
\end{equation}

Property (\ref{e3}) allows to apply the quasiclassical methods
to equation (\ref{e2}) with $\epsilon$ being the parameter of
the quasiclassical expansion.

Since the mean number of particles in the state $\hat\Phi$
is equal to \cite{B,LL}
$$N=\frac{(\hat\Phi,\hat\psi_l^+\hat\psi_l^-\hat\Phi)}
{(\hat\Phi,\hat\Phi)} = \frac{1}{\epsilon}\frac
{(\hat\Phi,\hat\phi_l^+\hat\phi_l^-\hat\Phi)}{(\hat\Phi,\hat\Phi)},
l=\overline{1,\infty}$$
and the quantity $(\hat\Phi,\hat\phi_l^+\hat\phi_l^-\hat\Phi)
/(\hat\Phi,\hat\Phi)$ is of order $\epsilon^0$ the quasiclassical
methods allow to construct approximate solutions of the equation
(\ref{e1}) at $\epsilon\to 0, N\to\infty,\epsilon N\to\alpha=\const$.

To the equation (\ref{e2}) there corresponds a classical infinite
dimesional Hamiltonian system with the Hamiltonian $H \left(
\frac{Q_j-iP_j}{\sqrt{2}},\frac{Q_j+iP_j}{\sqrt{2}}\right)$.

The results relating to statistical physics may be rigorously
justified. For a part of them such justification was given in
\cite{MS2}.

An analogous parameter arises in quantum field theory. Considering
formal asymptotic expansions over this small parameter does
not allow however to justify this asymptotics since the quantum
field theory itself does not have rigorous mathematical meaning.
Some approaches to this problem have been developed only in
 partial cases, see \cite{Si,GJ}. And without its solution
a justification of heuristic asymptotics is of course impossible.
Only the postulated perturbation theory series makes sense
in the quantum field theory up to now. One can expect that
a geometrical quantization coinciding with the "classical"
equations in the case when  all the commutators vanish
should be postulated as well.

The concept of geometrical quantization has been essentially
introduced in the works by Bohr, Sommerfeld, and de Broglie,
even before the works by Schr\" odinger,
Heisenberg, and Dirac. It is also a rapidly developping
concept in modern mathematics \cite{GS,KM}.

This paper is organized as follows. Section 2 contains definition of a
canonical operator in Fock space. We define objects of geometric
quantization -- Lagrangian manifolds with complex germs, and assign
elements of Fock sapace to these manifolds, justify a connection between
our and traditional definition of a canonical operator. These definitions
differ because of divergences in infinite dimensional case.
Traditional definition is not applicable in this case. In section 3
we define canonical transformation of a Lagrangian manifold with
complex germ and justify germ axioms. We show that a canonical operator
approximately satisfyes corresponding secondary - quantized equations.
Theorem is formulated in section 4. Section 5 contains the construction
of another asymptotics with the help of complex germ creation and
annihilation operators. Many examples are cited in sections 4, 5.
In section 6 we prove the theorem.

\section{A canonical operator}

We consider in the present paper finite dimensional isotropic
manifolds with infinite dimensional complex germ in
an infinite dimensional phase space (see the definition
below). We assign a canonical operator to such a manifold.
For the sake of simplicity we give a not absolutely invariant
definition of the canonical operator. This definition is
however sufficient for solving the Cauchy problem.

This canonical operator approximately satisfies the corresponding secondary
quantized equations provided that the initial conditions
for these equations correspond to isotropic finite dimensional
manifolds (or, in the stationary case, there exist
stable isotropic manifolds).

In specific examples such
isotropic manifolds usually have dimension $1$ or $0$.

\subsection{Objects of geometrical quantization}

Define now the objects of geometrical quantization.

By $\H^n$ denote the space of complex square summable symmetric functions
of $n$ variables $i_1,\dots,i_n\in \N$. Introduce in
$\H^n$ a scalar product by the following formula

$$(f,g) = \sum\limits_{i_1,\dots,i_n=1}^{\infty}
f_{i_1\dots i_n}^*g_{i_1\dots i_n}; \quad f,g\in\H^n.$$

Denote by $\H$ the Fock space $\oplus_{n=0}^{\infty}\H^n$,
and by $\hat\Phi^{(k)}\in\H^k$ the $k$-th component of $\hat\Phi\in\H$.

Consider creation and annihilation operators acting in $\H$ in the following
way:
\begin{eqnarray*}
(\hat\psi_j^-\hat\Phi)_{j_1\dots j_{k-1}}^{(k-1)} & = &
k^{1/2}\hat\Phi_{j_1\dots j_{k-1}j}^{(k)}, \\
(\hat\psi_j^+\hat\Phi)_{j_1\dots j_{k}}^{(k)} & = &
k^{-1/2}\sum\limits_{i=1}^{k}\hat
\Phi_{j_1\dots j_{i-1}j_{i+1}\dots j_k}^{(k-1)}\delta_{jj_i}, \\
& & j,j_l\in\N.
\end{eqnarray*}

By $\hat\Phi_0$ denote the following element of the space $\H$:
$\hat\Phi_0^{(0)} =1,
\hat\Phi_0^{(i)} =0, i\ge 1$.
As usual we call elements of the space $\H$ {\it state vectors
in the Fock presentation}.

We will also use in the present paper the presentation of elements in $\H$ as
the Berezin generating functionals \cite{B}.

To each element $\hat\Phi\in\H$ we assign the generating functional
\begin{eqnarray*}
\Phi(a_1^*,a_2^*,\dots) & = & \sum\limits_{n=0}^{\infty}
\frac{1}{\sqrt{n!}}\sum\limits_{i_1,\dots,i_n = 1}^{\infty}
\hat\Phi_{i_1\dots i_n}^{(n)}a_{i_1}^*\dots a_{i_n}^*, \\
& & a_j^* \in\C,\sum\limits_{j=1}^{\infty} \vert a_j^*\vert^2 <
\infty.
\end{eqnarray*}

In this presentation the creation operators $\hat\psi_j^+$
are multiplications by $a_j^*$, and the annihilation operators
$\hat\psi_j^-$ are derivation operators $\partial/\partial a_j^*$
\cite{B}.

We will also make use of the Schr\" odinger coordinate presentation
(or the $Q$-presentation). In this presentation we assign to each
element $\hat\Phi\in\H$ the functional
\begin{eqnarray*}
\Phi_Q(q_1,q_2,\dots) & = & \sum\limits_{n=0}^{\infty}\frac{1}
{\sqrt{n!}}\sum\limits_{i_1,\dots,i_n = 1}^{\infty}
\hat\Phi_{i_1\dots i_n}^{(n)} \\
& & \cdot\prod\limits_{l=1}^{n} [(q_{i_l}-\epsilon\partial/\partial
q_{i_l})(2\epsilon)^{-1/2}] \exp(-\sum\limits_{i=1}^{\infty}
q_i^2/2\epsilon).
\end{eqnarray*}

In this presentation the creation and annihilation
 operators $\hat\psi_j^{\pm}$
have the form $\hat\psi_j^{\pm}=
(q_j\mp\epsilon\partial/\partial q_j)/\sqrt{2
\epsilon}$.

By $\L$ denote the set

$$\L =\{(P_1,Q_1;P_2,Q_2;\dots):P_i\in \R, Q_i\in\R,
\sum\limits_{i=1}^{\infty}(P_i^2+Q_i^2)<\infty\}.$$

Let us define now finite dimensional isotropic manifolds and
corresponding Lagrangian manifolds with complex germs.

Let $\Lambda^k$ be a $k$-dimensional surface in $\L$,
and $\tau_1,\tau_2,\dots,\tau_k$ local coordinates on $\Lambda^k$.
For the sake of brevity denote the sequence $(P_1,P_2,P_3,\dots)$
by $P$, the sequence $(Q_1,Q_2,Q_3,\dots)$ by $Q$, and the
sequence $(\tau_1,\tau_2,\dots,\tau_k)$ by $\tau$.

\begin{definition}\label{d1}

A manifold $\Lambda^k=\{P=P(\tau),
Q=Q(\tau)\}$ is called {\em isotropic} if the following axioms hold
\begin{description}

\item[m1)]\label{am1} For any $\lambda_1,\dots,\lambda_k,
\lambda_i \in \{0,1,2,\dots\},i=\overline{1,k}$
the derivatives

$$\frac{\partial^{\lambda_1+\dots+\lambda_k}}{\partial\tau_1^{\lambda_1}
\dots\partial\tau_k^{\lambda_k}}P_j,
\frac{\partial^{\lambda_1+\dots+\lambda_k}}{\partial\tau_1^{\lambda_1}
\dots\partial\tau_k^{\lambda_k}}Q_j,\quad j=1,2,3,\dots $$
exist and the series
$$\sum\limits_{j=1}^{\infty}\left[\left(
\frac{\partial^{\lambda_1+\dots+\lambda_k}}{\partial\tau_1^{\lambda_1}
\dots\partial\tau_k^{\lambda_k}}P_j\right)^2 + \left(
\frac{\partial^{\lambda_1+\dots+\lambda_k}}{\partial\tau_1^{\lambda_1}
\dots\partial\tau_k^{\lambda_k}}Q_j\right)^2\right]$$
converges.

\item[m2)]\label{am2} If
$$\sum\limits_{j=1}^{\infty}\left[
\left(\sum\limits_{m=1}^k
\frac{\partial P_j}{\partial\tau_m}\xi_m\right)^2
+\left(\sum\limits_{m=1}^k
\frac{\partial Q_j}{\partial\tau_m}\xi_m\right)^2
\right] = 0,$$

for a $\xi_m\in\R$ then $\xi_m=0, m = \overline{1,k}$.

\item[m3)]\label{am3}

$$\sum\limits_{j=1}^{\infty}\left(
\frac{\partial P_j}{\partial\tau_m}
\frac{\partial Q_j}{\partial\tau_n}
- \frac{\partial P_j}{\partial\tau_n}
\frac{\partial Q_j}{\partial\tau_m}
\right) = 0, m,n=\overline{1,k}.$$
\end{description}
\end{definition}

Denote by $M^k$ the universal covering space of the isotropic manifold
$\Lambda^k$.

\begin{definition}

{\em A complex germ} is a set of planes $r(\tau),\tau\in M^k$
in the abstract space $\R^{\infty}\times\R^{\infty}$.
Each plane consists of the vectors

$$ \left(\begin{array}{c}
w_1(\alpha,\tau) \\ w_2(\alpha,\tau) \\ \vdots \\
z_1(\alpha,\tau) \\ z_2(\alpha,\tau) \\ \vdots
\end{array}\right),
\begin{array}{rcl}
w_i(\alpha,\tau) & = & \sum\limits_{j=1}^{\infty}
B_{ij}(\tau)\alpha_{j}, \quad i = \overline{1,\infty}, \\
z_i(\alpha,\tau) & = & \sum\limits_{j=1}^{\infty}
C_{ij}(\tau)\alpha_{j}, \quad i = \overline{1,\infty},\end{array}$$
where $\alpha$ ranges over infinite sequences
$\alpha = (\alpha_1,\alpha_2,\dots)
(\alpha_i \in \C,i=\overline{1,\infty},\sum\limits_{i=1}^{\infty}
\vert\alpha_i\vert^2<\infty)$ and
the following axioms hold:
\begin{description}

\item[r1)]\label{ar1}
Let $\tau^{'}, \tau^{''}$ be two points of the universal covering space
$M^k$ such that they project into one and the same point of $\Lambda^k$.
Then there exists an unitary operator $A(\tau^{'},\tau^{''})$ such that
$$B(\tau^{''})=B(\tau^{'})A(\tau^{'},\tau^{''}),
C(\tau^{''})=C(\tau^{'})A(\tau^{'},\tau^{''}),
 $$

\item[r2)]\label{ar2} For $a=\overline{1,k}, i=\overline{1,\infty}$
$$B_{ia}(\tau)=\partial P_i(\tau)/\partial\tau_a,
C_{ia}(\tau)=\partial Q_i(\tau)/\partial\tau_a.$$

\item[r3)]\label{ar3}
$$B^T(\tau)C(\tau) -
C^T(\tau)B(\tau) = 0. $$

(Here $B^T$ and $C^T$ are transpose matrices to $B$ and $C$
respectively.)

\item[r4)]\label{ar4}
$$C^+(\tau)B(\tau) -
B^+(\tau)C(\tau) = iL, $$
where $L$ is the diagonal matrix with first $k$ diagonal elements
equal to zero and all others to $1$.

\item[r5)]\label{ar5} For any $\lambda_1,\dots,\lambda_k;
\lambda_i\in\{0,1,2,\dots\},i=\overline{1,k}$
the derivatives
\begin{eqnarray*}
F_{mn}^{(\lambda_1,\dots,\lambda_k)} & = &
\frac{\partial^{\lambda_1+\dots +\lambda_k}}
{\partial\tau_1^{\lambda_1}\dots\partial\tau_k^{\lambda_k}}
\left(\frac{C_{mn}(\tau)+iB_{mn}(\tau)}{\sqrt{2}}\right), \\
G_{mn}^{(\lambda_1,\dots,\lambda_k)} & = &
\frac{\partial^{\lambda_1+\dots +\lambda_k}}
{\partial\tau_1^{\lambda_1}\dots\partial\tau_k^{\lambda_k}}
\left(\frac{C_{mn}(\tau)-iB_{mn}(\tau)}{\sqrt{2}}\right),
\end{eqnarray*}
exist, operators $G^{(\lambda_1,\dots,\lambda_k)}$ are bounded,
and operators $F^{(\lambda_1,\dots,\lambda_k)}$ are Hilbert-Schmidt
operators.

\item[r6)]\label{ar6}  The operator $(C(\tau)-iB(\tau))/\sqrt{2}$ has a
bounded inverse operator.

\end{description}
\end{definition}
A pair consisting of an isotropic manifold and a complex germ will be called
{\it a Lagrangian manifold with complex germ} since they form together
a germ of an infinite dimensional complex Lagrangian manifold.

\begin{remark}

 Let $r(\tau^{'}), r(\tau^{''})$ be the planes in
two points $\tau^{'},\tau^{''}$ of the universal covering space
$M^k$ that project into one and the same point of
the isotropic manifold $\Lambda^k$.
It follows then from the axiom r1 that these two planes differ
from each other only by a parametrization.

\end{remark}

\begin{remark}

In the finite dimensional case axioms r2-r4 are equivalent
to traditional axioms of the complex germ \cite{M1,M2,BD}.

\end{remark}

\begin{remark}

For a zero dimensional isortopic manifold the axioms of the
complex germ mean that the matrix
$$\left( \begin{array}{cc}
G & \bar F \\
F & \bar G
\end{array}\right) $$
is the matrix of a Berezin proper canonical transformation
(see  \cite{B}).

\end{remark}

Assign the function $\phi_l(\tau)=(Q_l(\tau)+iP_l(\tau))
/\sqrt{2}$ to each isotropic manifold $\Lambda^k
= \{P=P(\tau),Q=Q(\tau)\}$.

\subsection{Heuristic motivation for definition of a canonical
operator}
Consider a heuristic method to derive state vectors, which are approximate
solutions to secondary-quantized equations (\ref{e2}) and corresponds to
Lagrangian manifolds with complex germs.

Equation (\ref{e2}) in $Q$-presentation has the form of infinite-dimensional
Schr\" odinger equation,
\begin{equation}\label{A1}
i\epsilon\frac{\partial \Phi_{Q}}{\partial t}=
H\left(\frac{\hat Q_{j} - i\hat P_{j}}{\sqrt{2}},
\frac{\hat Q_{j} + i\hat P_{j}}{\sqrt{2}}\right)\Phi_{Q}
\end{equation}
where operators $\hat Q_{j}$ are multiplications by $q_{j}$,operators
$\hat P_{j}$ are derivation operators $-i\epsilon \partial/\partial q_{j}$.
Thus, $\hat Q_{j}$ and $\hat P_{j}$ are infinite-dimensional analogs
of coordinate and momentum operators, while $\epsilon$ is the analog of
Planck constant. As it was mentioned above, one can apply semiclassical
technique to equation (\ref{A1}) as $\epsilon$ tends to zero.

Semiclassical approximate solutions of the following type
\begin{equation}\label{A2}
\phi(q)\exp\left(\frac{i}{\epsilon} S(q)\right) ,
\end{equation}
where the sequence $(q_1,q_2,...)$ is denoted by $q$, $S$ is a real
functional, are widely used in physics ( see, for example,\cite{LL}).
Functionals (\ref{A2}) are of order $O(1)$ as $\epsilon \rightarrow 0$
at all $q$.

In this paper we consider yet another type of asymptotic solutions
to quantized equations. These solutions are not small and have rapidly
oscillating form only if the distance between the point $q$ and surface
\begin{equation}\label{A2_5}
 \{ Q(\tau),\tau \in \Lambda^k \}
\end{equation}
is of order $O(\epsilon^{1/2})$ as $\epsilon$ tends to zero,i.e., the
following quantity
\begin{equation}\label{A3}
\min_{\tau \in \Lambda^k} \sum_{i=1}^{\infty} (q_i-Q_i(\tau))^2
\end{equation}
is of order $O(\epsilon)$.If quantity (\ref{A3}) is of order
$O(\epsilon^{1-\delta}),\delta > 0,$
then these solutions are exponentially small.

First of all, consider the case of zero-dimensional Lagrangian manifold,
i.e. point
\newline
$(P_1,Q_1,P_2,Q_2,...) \in \L $ . Complex germ method leads to the
asymptotics of the following type
\begin{equation}\label{A4}
\Phi_Q(q) = {\cal F}\left(\frac{q-Q}{\sqrt{\epsilon}}\right)
\exp\left(\frac{i}{\epsilon}\sum_{i=1}^{\infty} P_i(q_i-Q_i)\right)
\end{equation}
Function $\cal F$ in formula (\ref{A4}) rapidly decays as its argument
tends to infinity and can be expressed through the complex germ
(see subsection 2.5 for more details).

Give Fock presentation for functional (\ref{A4}). First, consider for the
sake of simplicity case $P=Q=0$.If
\begin{equation}\label{A5}
 {\cal F}\left(\frac{q}{\sqrt{\epsilon}}\right)
= {\cal F}_{0}\left(\frac{q}{\sqrt{\epsilon}}\right) =
\exp(-\sum_{i=1}^{\infty} q_i^2/2\epsilon)
\end{equation}
then functional (\ref{A4}) corresponds to the vacuum state vector
$\hat\Phi_0$ in Fock presentation. If functional $\cal F$ is arbitrary,
then we can extract factor (\ref{A5}) from functional (\ref{A4}) and
present corresponding state vector in the form
\begin{equation}\label{A6}
 {\cal F}_{1}\left(\frac{\hat Q}{\sqrt{\epsilon}}\right)   \hat\Phi_0
\end{equation}
where ${\cal F}_{1}={\cal F}/{\cal F}_{0}$.
Operator ${\hat Q}/\sqrt{\epsilon}$ can be
expressed as a function of creation and annihilation operators
\begin{equation}\label{A7}
 \frac{\hat Q_i}{\sqrt{\epsilon}} =
 \frac{\hat \psi_i^+ + \hat \psi_i^-}{\sqrt{2}}
\end{equation}
which does not depend on $\epsilon$. Any function of operator (\ref{A7})
transform vacuum vector to $\epsilon$-independent element of the Fock
space, which corresponds to functional (\ref{A4}) if $P=Q=0$.

Consider now the general case, which can be reduced to the considered case.
Namely, functional (\ref{A4}) can be presented in the form
\begin{equation}\label{A8}
 \exp\left(\frac{i}{\epsilon}P_l q_l\right)
\exp\left(-Q_l \frac{\delta}{\delta q_l}\right)
  {\cal F}_{0}\left(\frac{q}{\sqrt{\epsilon}}\right)
\end{equation}
where $l=\overline{1,\infty}$. Expression (\ref{A8}) can be simplified
by obtaining Baker-Hausdorff formula
\begin{equation}\label{A9}
e^{A_1+A_2}=e^{A_1}e^{A_2}e^{-[A_1,A_2]/2}
\end{equation}
for operators
$$
A_1=iP_l q_l/\epsilon,A_2=-Q_l \delta/\delta q_l, l=\overline{1,\infty}
$$
{}From
$$
A_1+A_2 = \frac{1}{\sqrt{\epsilon}} \sum_{l=1}^{\infty}
(\hat\psi_l^+ \phi_l - \hat\psi_l^- \phi_l^* )
$$
we obtain the following state vector in Fock presentation, corresponding to
functional (\ref{A8}) :
\begin{equation}\label{A10}
U_{\phi} \exp(-\frac{i}{2\epsilon}\sum_{l=1}^{\infty}P_l Q_l) \hat Y
\end{equation}
where
$$
U_{\phi} = \exp( \frac{1}{\sqrt{\epsilon}} \sum_{l=1}^{\infty}
(\hat\psi_l^+ \phi_l - \hat\psi_l^- \phi_l^* ))
$$
$\hat Y$ is $\epsilon$-independent state vector, which can
be expressed through creation operators \cite{B}:
$$
\hat Y=Y(\hat \psi^+)\hat \Phi_0
$$
$$
Y(\hat \psi^+)=\sum_{n=0}^{\infty}\frac{1}{\sqrt{n!}} Y^{(n)}_{i_1,...,i_n}
 \hat \psi^+_{i_1}...\psi^+_{i_n}
$$
{}From Baker-Hausdorff formula and the following commutation relations
\cite{B}
\begin{equation}\label{A11}
U_{\phi}^{-1} \hat \psi_l^+ U_{\phi} =
\hat \psi_l^+ + \phi_l^*/\sqrt{\epsilon},\quad
U_{\phi}^{-1} \hat \psi_l^- U_{\phi} =
\hat \psi_l^- + \phi_l/\sqrt{\epsilon}
\end{equation}
which can be verified by using commutation relations between creation
and annihilation operators, we obtain that vector (\ref{A10}) can be
expressed in a form
\begin{equation}\label{A12}
\Psi_{\phi,Y}=Y(\hat \psi^+ -\phi^*/\sqrt{\epsilon})
\exp\left(\frac{1}{\epsilon} [ g + \sum_{l=1}^{\infty} \phi_l
(\sqrt{\epsilon}\hat\psi_l^+ - \phi_l^* )]\right) \hat\Phi_0
\end{equation}
where
\begin{equation}\label{A12_1}
g=\phi_l^* \phi_l /2 + ( \phi_l^* \phi_l^* - \phi_l \phi_l )/4,
l=\overline{1,\infty}
\end{equation}
We shall usually use Gaussian functionals $Y$
\begin{equation}\label{A12_2}
Y(\hat\psi^+)=c \exp(\hat \psi_i^+ M_{ij} \hat \psi_j^+ /2),
i,j=\overline{1,\infty}
\end{equation}
where $c$ is a constant, $M$ is a Hilbert-Schmidt operator,$\|M\| < 1$.

Thus, we have obtained state vector (\ref{A12}), corresponding to
zero-dimensional
 isotropic manifold. We shall assign matrix $M$ to each complex
germ in subsection 2.3. Discussion about connection between complex
germ theories in $Q$-presentation and in Fock space in more details
will be presented in subsection 2.5.

Consider now asymptotics in $Q$-presentation, which are peaked in the
vicinity of the surface (\ref{A2_5}) if $k>0$. In this subsection we
discuss partial case of manifold $\Lambda^k$, which satisfies the
conditions
\begin{equation}\label{A13}
Q_{k+1}=Q_{k+2}=...=P_{k+1}=P_{k+2}=...=0
\end{equation}
General case will be considered in subsection 2.5.

Complex germ asymptotics in a case (\ref{A13}) has the form
\begin{equation}\label{A14}
{\cal F}\left(q_1,...,q_k;\frac{q_{k+1}}{\sqrt{\epsilon}},
\frac{q_{k+2}}{\sqrt{\epsilon}},...\right)
\exp\left(\frac{i}{\epsilon}S(q_1,...,q_k)\right),
\end{equation}
where $\cal F$ rapidly decays at infinity, $S$ is a real function.
Isotropic manifold,corresponding to asymptotics (\ref{A14}), can be
parametrized by coordinates $\tau_a=Q_a,a=\overline{1,k}$
$$
Q_a(\tau)=\tau_a, P_a(\tau)=\partial S/\partial \tau_a
$$
In order to give Fock presentation of functional (\ref{A14}) it is
convenient to reduce it to functionals (\ref{A4}). Namely, consider
the following superposition of expressions (\ref{A13})
$$
\int \frac{d\tau}{\epsilon^{k/2}} f\left(\tau_1,...,\tau_k;
\frac{q_1-\tau_1}{\sqrt{\epsilon}},...,
\frac{q_k-\tau_k}{\sqrt{\epsilon}};
\frac{q_{k+1}}{\sqrt{\epsilon}},
\frac{q_{k+2}}{\sqrt{\epsilon}},...\right)
$$
\begin{equation}\label{A15}
\times\exp\left(\frac{i}{\epsilon}\left[S(\tau)+\sum_{a=1}^{k}
\frac{\partial S}{\partial \tau_a}(\tau)(q_a-\tau_a)\right]\right)
\end{equation}
and choose rapidly decaying at infinity function $f$ in order to make
functional (\ref{A15}) approximately equal to (\ref{A14}). Consider the
following substitution in integral (\ref{A15}):
$$
\xi_a=(q_a-\tau_a)/\sqrt{\epsilon}, a=\overline{1,k},
$$
which transforms it to the following expression:
$$
\int d\xi_1...d\xi_k
f\left(q_1-\xi_1\sqrt{\epsilon},...,q_k-\xi_k\sqrt{\epsilon};
\xi_1,...,\xi_k;
\frac{q_{k+1}}{\sqrt{\epsilon}},
\frac{q_{k+2}}{\sqrt{\epsilon}},...\right)
$$
\begin{equation}\label{A16}
\times\exp\left(\frac{i}{\epsilon}\left[S(q-\xi\sqrt{\epsilon})
+\sum_{a=1}^{k}
\frac{\partial S}{\partial q_a}(q-\xi\sqrt{\epsilon})
\xi_a\sqrt{\epsilon}\right]\right)
\end{equation}
If $\epsilon\rightarrow 0 , \xi=\const$ then the exponent in formula
(\ref{A16}) has the form
$$
\exp\left[-\frac{i}{2}\xi_a
\frac{\partial^2 S(q)}{\partial q_a\partial q_b}\xi_b + \frac{i}{\epsilon}
S(q)\right]
$$
Taking into account that $f$ rapidly decays, we obtain that formulas
(\ref{A15}) and (\ref{A14}) are approximately equal if
$$
{\cal F}\left(q_1,...,q_k;
\frac{q_{k+1}}{\sqrt{\epsilon}},
\frac{q_{k+2}}{\sqrt{\epsilon}},...\right)
$$
\begin{equation}\label{A17}
= \int d^k\xi
f\left(q_1,...,q_k;
\xi_1,...,\xi_k;
\frac{q_{k+1}}{\sqrt{\epsilon}},
\frac{q_{k+2}}{\sqrt{\epsilon}},...\right)
\exp\left[-\frac{i}{2}\xi_a
\frac{\partial^2 S(q)}{\partial q_a\partial q_b}\xi_b \right]
\end{equation}
Of course, there are many such functions $f$.

Functional (\ref{A15}) has the following form in Fock presentation
\begin{equation}\label{A18}
\int \frac{d\tau}{\epsilon^{k/2}}
Y(\hat \psi^+ -\phi^*/\sqrt{\epsilon})
\exp\left(\frac{1}{\epsilon} [ g(\tau) + iS(\tau)
+ \sum_{l=1}^{\infty} \phi_l(\tau)
(\sqrt{\epsilon}\hat\psi_l^+ - \phi_l^*(\tau) )]\right) \hat\Phi_0
\end{equation}

In formula (\ref{A18}) $g(\tau)$ has the form (\ref{A12_1}). We shall
usually use functionals $Y$ of the form (\ref{A12_2}). Vector (\ref{A18})
will be multiplied by $\epsilon^{k/4}$ for making its norm of order $O(1)$.
Notice also that
$$
g(\tau) + iS(\tau) = \int_{\tau^{(0)}}^{\tau} \phi_l d\phi_l^*  +
g(\tau^{(0)}) + iS(\tau^{(0)})
$$
We shall define a canonical operator by formula, analogous to
(\ref{A18}) in subsection 2.4; some auxiliary lemmas will be proved
in subsection 2.3. Connection between the  expression (\ref{A18})
and traditional canonical operator, corresponding to arbitrary Lagrangian
manifold with complex germ,
will be discussed in subsection 2.5.
\subsection{Some auxiliary lemmas}

We are going now to prove some lemmas.

\begin{lemma}\label{l1}

Let $L$ be an operator in $l^2$ that maps a sequence
$(\xi_1,\xi_2,\xi_3,\dots)$ into the sequence
$(0,0,\dots,0,\xi_{k+1},\xi_{k+2},\dots)$, and let $Y$
be a bounded operator in $l^2$ satisfying the following
properties:
\begin{enumerate}
\item\label{l1ia}
$(\xi,Y\xi)\ge 0,\xi\in l^2;$
\item\label{l1ib}
if $L\xi=0,(\xi,Y\xi)=0,$  then $\xi=0$.
\end{enumerate}
Then there exists a positive $\kappa$, such that
$$(\xi,(L+Y)\xi)\ge\kappa(\xi,\xi)$$
for any $\xi$.

\end{lemma}

\sbp{Proof.}

It follows from the assumptions of the lemma that

\begin{equation}\label{e4}
(\xi,(L+Y)\xi)\ge(L\xi,L\xi).
\end{equation}

On the other hand, the property \ref{l1ib} implies that
for $L\xi = 0$
$$(\xi,Y\xi)\ge\sigma(\xi,\xi),\sigma > 0.$$

It follows then that
\begin{eqnarray}\label{e5}
(\xi,(L+Y)\xi) & = & ((E-L)\xi,Y(E-L)\xi)+ (L\xi,Y(E-L)\xi) \nonumber \\
& & +((E-L)\xi,YL\xi)+(L\xi,(YL+L)\xi) \nonumber \\
& \ge & \sigma(\xi,\xi)-2\|Y\|\sqrt{(L\xi,L\xi)(\xi,\xi)}
-(\|Y\|+\sigma)(L\xi,L\xi),
\end{eqnarray}
where $\| Y\| = \sup\limits_{\|\xi\|=1}\|Y\xi\|.$

Suppose that for any $\delta > 0$ there exists a vector $\xi\in l^2$
such that $(\xi,(L+Y)\xi)<\delta(\xi,\xi)$. It follows then from
(\ref{e4}) that $(L\xi,L\xi)<\delta(\xi,\xi)$, and it follows from
(\ref{e5}) that
$$ (\xi,(L+Y)\xi)\ge(\sigma-(\|Y\|+\sigma)\delta -2\|Y\|
\sqrt{\delta})(\xi,\xi).$$
This inequality for $\delta$ small enough contradicts the
condition $(\xi,(L+Y)\xi)<\delta(\xi,\xi)$. The contradiction obtained
proves Lemma \ref{l1}.

Denote by $W_{ab}(\tau),\tau\in\Lambda^k, a,b=\overline{1,k}$
the inverse matrix to the matrix $\sum\limits_{i=1}^{\infty}
\frac{\partial\phi_i^*}{\partial\tau_a}\frac{\partial\phi_i}
{\partial\tau_b}$, which is invertable by axiom
m2. We set
\begin{equation}\label{e7}
M_{ij}(\tau)=((C+iB)(\tau)(C-iB)^{-1}(\tau))_{ij}
-\sum\limits_{a,b=1}^k\frac{\partial\phi_i}{\partial\tau_a}
W_{ab}(\tau)\frac{\partial\phi_j}{\partial\tau_b}
\end{equation}

\begin{lemma}\label{l2}

The operator $M$ has the following properties
\begin{enumerate}
\item $M$ is a Hilbert-Schmidt operator;
\item $\|M\|<1$;
\item $\sum\limits_{j=1}^{\infty}M_{ij}\partial\phi_j^*
/\partial\tau_a = 0,\quad a=\overline{1,k}$.
\end{enumerate}
\end{lemma}

\sbp{Proof.} Since the operators $C+iB$ and $\sum_{a,b=1}^{k}
\frac{\partial\phi_i}{\partial\tau_a}W_{ab}\frac{\partial\phi_j}
{\partial\tau_b}$ are Hilbert-Schmidt operators, and the operator
$(C-iB)^{-1}$ is a bounded operator, we have that $M$ is a
Hilbert-Schmidt operator.

Axiom r2 and definition (\ref{e7}) of $M$ imply that
$$\sum\limits_{j=1}^{\infty}M_{ij}\partial\phi_j^*/
\partial\tau_a = 0, \quad a=\overline{1,k}.$$

Expression (\ref{e7}) imply as well that
\begin{eqnarray*}
E-M^+M & = & ((C-iB)^+)^{-1}\left((C+iB)^+
\sum\limits_{a,b=1}^k\sum\limits_{l=1}^{\infty}
\frac{\partial\phi_l}{\partial{\tau_a}}\right. \\
& & \left.\cdot W_{ab}\frac{\partial\phi_l}{\partial\tau_b}
(C+iB)+2L\right)(C-iB)^{-1}
\end{eqnarray*}

Lemma \ref{l1} implies that for some $\rho > 0$
$(\xi,(E-M^{+}M)\xi)>\rho(\xi,\xi)$. It follows then that
$\|M\|<1$. Lemma \ref{l2} is proved.

\subsection{Definition of a canonical operator}

Consider an aggregate consisting of a number $\epsilon>0$,
a sequence $\phi_j, \sum\limits_{j=1}^{\infty}\vert\phi_j\vert^2
<\infty$, and a Hilbert-Schmidt operator $M:l^2\to l^2,
\|M\|<1$. Assign the following element of $\H$ to such an
aggregate:
\begin{eqnarray}\label{e8}
\hat\Phi_{\phi,M} & = & \exp\{\frac{1}{\epsilon}
\phi_j(\sqrt{\epsilon}\hat\psi_j^+-\phi_j^{*})\nonumber \\
& & +\frac{1}{2\epsilon}(\hat\psi_i^+\sqrt{\epsilon} -\phi_i^*)
M_{ij}
(\hat\psi_j^+\sqrt{\epsilon} -\phi_j^*) \}\hat\Phi_0,
\end{eqnarray}
where we mean summing over repeating indices, $i,j = \overline{1,\infty}$.

\begin{remark}

Expanding the exponent in formula (\ref{e8}) into the power series
over the variable $\epsilon$ we obtain the $l$-th component
of $\hat\Phi_{\phi,M}\in\H$ in the following form:
\begin{eqnarray}\label{e9}
(\hat\Phi_{\phi,M})_{i_1\dots i_l}^{(l)} & = &
\sum\limits_{k=0}^{[l/2]} \frac{c}{\sqrt{l!}2^kk!}
a_{i_1}\dots a_{i_l} \nonumber \\
& & \cdot\sum\limits_{1\le j_1\ne\dots\ne j_{2k}\le l}
\frac{M_{i_{j_1}i_{j_2}}\dots M_{i_{j_{2k-1}}i_{j_{2k}}}}
{a_{i_{j_1}}\dots a_{i_{j_{2k}}}}
\end{eqnarray}
where
\begin{eqnarray*}
a_m & = & (\phi_m-M_{mn}\phi_n^*)/\sqrt{\epsilon}, \quad
m,n = \overline{1,\infty} \\
c & = & \exp\{-\frac{1}{\epsilon}\phi_j^*\phi_j
- \frac{1}{2\epsilon}\phi_i^*M_{ij}\phi_j^*\}, \quad
i,j = \overline{1,\infty}
\end{eqnarray*}
\end{remark}

Consider a Lagrangian manifold with complex germ satisfying
the quantization conditions
\begin{equation}\label{e10}
\frac{1}{2\pi \epsilon}\oint\limits_lP_idQ_i = n_{l},
\quad n_{l}\in{\Z},i=\overline{1,\infty}
\end{equation}
$l$ is an arbitrary closed path on the isotropic manifold $\Lambda^k$.

Let $f$ be an infinitely differentiable function
with compact support on the
isotropic manifold $\Lambda^k$, and let $\tau^{(0)}$ be an arbitrary
point on $\Lambda^k$.

Denote
$$ g= \phi_l^*(\tau^{(0)})\phi_l(\tau^{(0)})/2 +
(\phi_l^*(\tau^{(0)})\phi_l^*(\tau^{(0)})-
\phi_l(\tau^{(0)})\phi_l(\tau^{(0)}))/4,
l=\overline{1,\infty},$$
$$ F(\tau)=\frac{C(\tau)+iB(\tau)}{\sqrt{2}},
   G(\tau)=\frac{C(\tau)-iB(\tau)}{\sqrt{2}},
   \tau \in \Lambda^{k}. $$

We assign the following element of $\H$
to the Lagrangian manifold with complex germ $[\Lambda^k,r]$,
the function $f$, the point $\tau^{(0)}\in\Lambda^k$,
and the number $\epsilon>0$:
\begin{equation}\label{e11}
\begin{array}{rcl}
\hat\K_{[\Lambda^k,r],\tau^{(0)}}^{\epsilon} f & = &
 \int\limits_{\Lambda^k}
\frac{d\sigma f(\tau)}{(2\pi)^{k/2}\epsilon^{k/4}} \\
& & \cdot
\frac{\exp(\frac{1}{\epsilon}[g+\int_{\tau^{(0)}}^{\tau}
\phi_i(\tau^{'})d\phi_i^*(\tau^{'})])}
{\sqrt[4]{\det G^+(\tau)G(\tau)}}
\hat\Phi_{\phi(\tau),M(\tau)}^{\epsilon} \\
& & i =\overline{1,\infty};
\end{array}
\end{equation}
where $M$ is defined by formula (\ref{e7}),
$d\sigma$ is the following measure on $\Lambda^{k}$
\[
 d\sigma =      \sqrt{\det
\frac{\partial\phi_i^*}{\partial\tau_a}(\tau)
\frac{\partial\phi_i}{\partial\tau_b}(\tau)} d\tau,
\]
$d\tau \equiv d\tau_1\dots d\tau_k; a,b=\overline{1,k},
i=\overline{1,\infty}.$
This measure does not depend on the choice of local coordinates
$\tau_1,...,\tau_k$ on $\Lambda^k$.

\begin{remark}

Since the manifold $\Lambda^k$ is isotropic
the integral $\int_{\tau^{(0)}}^{\tau}\phi_i(\tau^{'})
d\phi_i^{*}(\tau^{'})$ does not depend locally on the path.

\end{remark}

\begin{remark}

The integrand in (\ref{e11}) is a one-valued function due to the
quantization conditions (\ref{e10}).

\end{remark}

\begin{remark}  It follows from the axiom r4 that
\begin{equation}\label{e6}
(C-iB)^+(C-iB)=2L+(C+iB)^+(C+iB)
\end{equation}
Formula (\ref{e6}) implies that the operator
$ G^+(\tau)G(\tau)-E $
is an operator of trace class, and therefore Fredholm
determinant for the operator
$ G^+(\tau)G(\tau) $
is defined.

\end{remark}

\begin{remark}

By the axiom r6
$ det[G^+(\tau)G(\tau)]>0 $

\end{remark}

\subsection{Connection with the traditional definition of a canonical
operator}

Consider now the relation between definition (\ref{e11})
of a canonical operator and the traditional definition
 of a canonical operator \cite{M1,BD}. Write down the vectors
(\ref{e8}) and (\ref{e11}) in $Q$-presentation. Consider an auxiliary
presentation with the wave function depending on
$l$ momenta and infinite number of coordinates. This
presentation will allow us to involve the case of focal
points into consideration.

Let $I$ be a finite set of positive integers,
$I=\{i_1,i_2,\dots,i_l\}$. Denote by $\Lambda$ the
diagonal matrix with $\Lambda_{jj}=1$ for $j\not\in I$,
$\Lambda_{jj} = -j$ for $j\in I$. We assign the following
presentation ($I$-presentation) of the vectors of $\H$
to each set $I$.

To each element of $\H$ we assign the functional
\begin{eqnarray*}
\Phi_I(q_1^I,q_2^I,\dots) & = &
\int\frac{dq_{i_1}\dots dq_{i_l}}{(2\pi\epsilon)^{l/2}}
\exp\left(-\frac{i}{\epsilon}\sum\limits_{j\in I}q_j^Iq_j
\right) \\
& & \cdot\Phi_Q(q_1^I,\dots,q_{i_1-1}^I,q_{i_1},q_{i_1+1}^I,
\dots,q_{i_l-1}^I,q_{i_l},q_{i_l+1}^I,\dots)
\end{eqnarray*}

The $Q$-presentation is a special case of the $I$-presentation for
$I=\emptyset$.

In order to give $I$-presentations of the vectors (\ref{e8}) and (\ref{e11})
use the formula \cite{B}
$$\Phi(a^*) = \int K_z(a^*)\Phi(z^*)e^{-z^*z}
\prod dz^*dz.$$
Here $K_z(a^{*})$ is the generating functional corresponding
to the vector $K_z = \exp(\sum\limits_{j=1}^{\infty}
z_j\hat\psi_{j}^{+})\hat\Phi_0$, the measure in the
functional integral is defined in \cite{B}, and $\hat\Phi$ is an
arbitrary element of $\H$.

Introduce following notation
$$ P_j^I(\tau) = - Q_j(\tau), j\in I, P_j^I(\tau) = P_j(\tau),
j \not\in I; $$
$$ Q_j^I(\tau) = P_j(\tau), j\in I, Q_j^I(\tau) = Q_j(\tau),
j \not\in I; $$
$$ \phi_j^I(\tau) = (Q_j^I(\tau)+iP_j^I(\tau))/\sqrt{2}
= \Lambda_{jm}\phi_m(\tau), j,m=\overline{1,\infty}, $$
$$z_j^I = \Lambda_{jm}z_m,z_j^{*I} = \Lambda_{jm}^* z_m^*,
j,m = \overline{1,\infty}.$$

The functional
\begin{eqnarray*}
(K_z)_I(q_1^I,q_2^I,\dots) & = &
\exp\left\{\frac{1}{\sqrt{2\epsilon}}z_m^I
\left(q_m^I-\epsilon\frac{\partial}{\partial q_m^I}\right)\right\}
\exp\left\{-\frac{1}{2\epsilon}q_m^Iq_m^I\right\} \\
& = & \exp\left\{-\frac{1}{2}z_m^Iz_m^I+\sqrt{\frac{2}{\epsilon}}
z_m^Iq_m^I-\frac{1}{2\epsilon}q_m^Iq_m^I\right\}, \\
& & m=\overline{1,\infty}
\end{eqnarray*}
corresponds to the vector $\hat K_z$ in the $I$-presentation.

This expression follows from the Baker-Hausdorff formula
$$e^{A_1+A_2} = e^{A_1}e^{A_2}e^{-\frac{1}{2}[A_1,A_2]},
A_1=\frac{1}{\sqrt{2\epsilon}}z_m^Iq_m^I,
A_2 = -\sqrt{\frac{\epsilon}{2}}z_m^I\frac{\partial}
{\partial q_m^I} $$
and the equality
$$\exp\left\{a_i\frac{\partial}{\partial q_i^I}\right\}
f(q_1^I,q_2^I,\dots) = f(q_1^I+a_1,q_2^I+a_2,\dots).$$
It now follows  that the functional corresponding to the vector
$\hat\Phi_{\phi,M}$ has the following form in $I$-presentation:
\begin{eqnarray*}
(\Phi_{\phi,M})_I(q_1^I,q_2^I,\dots) & = &
\int\prod dz^*dz \exp\left(-z_l^{*I}z_l^I -\frac{1}{2}
z_l^Iz_l^I+\sqrt{\frac{2}{\epsilon}}z_l^Iq_l^I
-\frac{1}{2\epsilon}q_l^Iq_l^I \right. \\
& & +\frac{1}{2\epsilon}(z_l^{*I}\sqrt{\epsilon}-\phi_l^{*I})
M_{lm}^I(z_m^{*I}\sqrt{\epsilon}-\phi_m^{*I}) \\
& & \left.+\frac{1}{\epsilon}\phi_l^I(\sqrt{\epsilon} z_l^{I*}
-\phi_l^{I*})\right), \quad l,m=\overline{1,\infty},
\end{eqnarray*}
where $M^I = \Lambda M \Lambda$.

By \cite{B} this integral equals to
\begin{equation}\label{e12}
\begin{array}{rcl}
(\Phi_{\phi,M})_I(q_1^I,q_2^I,\dots) & = &
\exp\left\{\frac{1}{\epsilon}\frac{\phi_l^I-\phi_l^{*I}}{\sqrt{2}}
\left(q_l^I-\frac{\phi_l^I+\phi_l^{*I}}{\sqrt{2}}\right)\right. \\
& & +\frac{i}{2\epsilon}
\left(q_m^I-\frac{\phi_m^I+\phi_m^{I*}}{\sqrt{2}} \right)
A_{mn}^{I}
\left(q_n^I-\frac{\phi_n^I+\phi_n^{I*}}{\sqrt{2}}\right)
- \frac{\phi_l^{*I}\phi_l^{I}}{2\epsilon} \\
& & \left. +\frac{\phi_l^I\phi_l^I-\phi_l^{I*}\phi_l^{I*}}
{4\epsilon}\right\}\frac{1}{\sqrt{\det (E+M^I)}} \\
& & l,m,n = \overline{1,\infty} \\
A^I & = & i(E-M^I)(E+M^I)^{-1}
\end{array}
\end{equation}

\begin{remark}

Not any element of $\H$ of form (\ref{e8}) may be presented in form
(\ref{e12}). Indeed, $\det (E+M^I)$ exists only if the operator
$M^I$ is an operator of trace class, and the (\ref{e8})
may be defined in other cases as well. For example,
$$\exp\left\{\frac{1}{4}\sum\limits_{n=1}^{\infty}
\frac{1}{n}(\hat\psi_n^+)^2\right\}
\hat\Phi_0\in\H$$
but this vector may not be written down in the $Q$-presentation.
It means that the definition of a canonical operator in
Fock presentation (\ref{e11})
is more general than the corresponding definition
in $Q$-presentation.

\end{remark}

Consider now the vector (\ref{e11}) in the $I$-presentation.

By (\ref{e11}) and (\ref{e12}) we have
$$(\K_{[\Lambda^k,r],\tau^{(0)}}^{\epsilon}f)_I
(q_1^I,q_2^I,\dots) = \int_{\Lambda^k}
\frac{d\tau f(\tau)\sqrt{\det
\frac{\partial\phi_l^*}{\partial\tau_a}(\tau)
\frac{\partial\phi_l}{\partial\tau_b}(\tau)}}
{(2\pi)^{k/2}\epsilon^{k/4}} $$

$$\cdot \frac{\exp\left(\frac{i}{\epsilon}
\int_{\tau^{(0)}}^{\tau}P_l^I(\tau^{'})dQ_l^I(\tau^{'})+
\frac{i}{\epsilon}P_l^I(\tau)(q_l^I-Q_l^I(\tau)) +
\frac{i}{2\epsilon}(q_l^I-Q_l^I(\tau))A_{lm}^I(\tau)
(q_m^I-Q_m^I(\tau))\right)}
{\sqrt[4]{\det\left[\left(\frac{C(\tau)-iB(\tau)}{\sqrt{2}}
\right)^+\left(\frac{C(\tau)-iB(\tau)}{\sqrt{2}}\right)\right]}
\sqrt{\det(E+M^I(\tau))}} $$

\begin{equation}\label{e13}
\cdot\exp\left(-\frac{i}{\epsilon}\sum\limits_{j\in I}
P_j(\tau^{(0)})Q_j(\tau^{(0)})\right),
\end{equation}
$$l,m = \overline{1,\infty}, a,b=\overline{1,k} $$
where $A_{lm}^{I}$ is the matrix of the operator $A^I(\tau)
=i(E+\Lambda M(\tau)\Lambda)^{-1}(E-\Lambda M(\tau)\Lambda)$,
and  $M(\tau)$ is determined by formula (\ref{e7}).

In particular, we have in the $Q$-representation
$$(\K_{[\Lambda^k,r],\tau^{(0)}}^{\epsilon}f)_Q
(q_1,q_2,\dots) = \int_{\Lambda^k}
\frac{d\tau f(\tau)\sqrt{\det
\frac{\partial\phi_l^*}{\partial\tau_a}(\tau)
\frac{\partial\phi_l}{\partial\tau_b}(\tau)}}
{(2\pi)^{k/2}\epsilon^{k/4}} $$

$$\cdot \frac{\exp\left(\frac{i}{\epsilon}
\int_{\tau^{(0)}}^{\tau}P_l(\tau^{'})dQ_l(\tau^{'})+
\frac{i}{\epsilon}P_l(\tau)(q_l-Q_l(\tau)) +
\frac{i}{2\epsilon}(q_l-Q_l(\tau))A_{lm}(\tau)
(q_m-Q_m(\tau))\right)}
{\sqrt[4]{\det\left[\left(\frac{C(\tau)-iB(\tau)}{\sqrt{2}}
\right)^+\left(\frac{C(\tau)-iB(\tau)}{\sqrt{2}}\right)\right]}
\sqrt{\det(E+M(\tau))}} $$

\begin{equation}\label{e14}
l,m = \overline{1,\infty}, a,b=\overline{1,k}
\end{equation}

$$A(\tau) = i(E+M(\tau))^{-1}(E-M(\tau)).$$

We present now a heuristic justification of a relation between
the definition of a canonical operator in $Q$-presentation
(\ref{e14}) and the traditional definition of a canonical operator
\cite{M1,BD}.

Cover the support of the function $f$ on the isotropic manifold
$\Lambda^k$ by a finite number of domains $\Omega_{\alpha},
\alpha=\overline{1,M}$ such that each domain has a one-to-one
projection on one of the coordinate planes of the form

$$ P_j = 0, j\not\in I_{\alpha},Q_j = 0, j\in I_{\alpha}$$
$$I_{\alpha} = \{i_1^{\alpha},i_2^{\alpha},\dots,
i_{l_{\alpha}}^{\alpha}\}\subset \N. $$

Consider a partition of unity $1=\sum\limits_{\alpha=1}^{M}
e_{\alpha}(\tau)$, where functions $e_{\alpha}(\tau)$ are
infinitely differentiable with a support in $\Omega_{\alpha}$.
Introduce a notation $f_{\alpha}(\tau)=f(\tau)e_{\alpha}(\tau)$.

The canonical operator in the $Q$-presentation has the form
$$
\left(\K_{[\Lambda^k,r],\tau^{(0)}}^{\epsilon}f\right)_Q
(q_1,q_2,\dots)  =
\sum_{\alpha=1}^M\int\frac{dq_{i_1}^{\alpha}\dots dq_{i_{l_{\alpha}}}
^{\alpha}}
{(2\pi\epsilon)^{l_{\alpha}/2}}
\cdot \exp\left(-\frac{i}{\epsilon}\sum\limits_{j\in I_{\alpha}}
q_j^I q_j\right)
$$
\begin{equation}\label{e15}
 \cdot\left(\K_{[\Lambda^k,r],\tau^{(0)}}^{\epsilon}f_{\alpha}
\right)_{I_{\alpha}}
\left(q_1^I,\dots,
q_{i_1^{\alpha}-1}^I,q_{i_1^{\alpha}},
q_{i_1^{\alpha}+1}^I,\dots,
q_{i_l^{\alpha}-1}^I,q_{i_l^{\alpha}},
q_{i_l^{\alpha}+1}^I,\dots\right)
\end{equation}

Calculate now the integral in  formula (\ref{e13}) that
expresses $(\K_{[\Lambda^{k},r],\tau^{(0)}}^{\epsilon}
f_{\alpha})_{I_{\alpha}}$.

Introduce following notation:
$$\begin{array}{rclrclc}
B_{mn}^I & = & -C_{mn}, & C_{mn}^I & = & B_{mn}, & m\in I \\
B_{mn}^I & = & B_{mn}, & C_{mn}^I & = & C_{mn}, & m\not\in I
\end{array}$$
and $S_{ab}$ is the inverse matrix to the matrix
$$\frac{\partial Q_i^I}{\partial \tau_a}A_{ij}^I
\frac{\partial Q_j^I}{\partial\tau_b}
- \frac{\partial P_i^I}{\partial \tau_a}
\frac{\partial Q_i^I}{\partial\tau_b},
\quad a,b,=\overline{1,k}, \quad i,j =\overline{1,\infty}.$$

The following statement holds.

\begin{statement}
\begin{enumerate}
\item
\begin{eqnarray}\label{e16}
(B^I(C^I)^{-1})_{mn} & = &
A_{mn}^I - \left(A_{mr}^I\frac{\partial Q_r^I}{\partial\tau_a}
-\frac{\partial P_r^I}{\partial\tau_a}\right)
 S_{ab}\left(\frac{\partial Q_j^I}{\partial\tau_b}
A_{jn}^I - \frac{\partial P_n^I}{\partial\tau_b}\right),
\nonumber \\ & &
j,m,n,r = \overline{1,\infty}, \quad a,b = \overline{1,k}
\end{eqnarray}

\item
\begin{eqnarray}\label{e17}
\frac{\det\left(\frac{\partial\phi_m^*}{\partial\tau_a}
\frac{\partial\phi_m}{\partial\tau_b}\right)}
{\det(E+M^I)} & = &
\det\left[\frac{1}{2}(E-iB^I(C^I)^{-1})\right] \nonumber \\
& & \cdot \det\left[\frac{1}{i}\left(
\frac{\partial Q_m^I}{\partial\tau_a}
A_{mn}^I
\frac{\partial Q_n^I}{\partial\tau_b} -
\frac{\partial P_m^I}{\partial\tau_a}
\frac{\partial Q_m^I}{\partial\tau_b}\right)\right] \\
& & m,n = \overline{1,\infty},\quad a,b =\overline{1,k} \nonumber
\end{eqnarray}
\end{enumerate}
\end{statement}

\sbp{Proof.}

Property (\ref{e16}) follows from the following expressions for
$B^I(C^I)^{-1}$ and $A^I$:
$$B^I(C^I)^{-1} = i(E+N^I)^{-1}(E-N^I),$$
$$N_{mn}^I = M_{mn}^I+\frac{\partial\phi_m^I}{\partial\tau_a}
W_{ab}\frac{\partial\phi_n^I}{\partial\tau_b},
\quad a,b =\overline{1,k}, \quad m,n = \overline{1,\infty},$$
where $W_{ab}$ is the inverse matrix to the matrix
$\frac{\partial\phi_m^*}{\partial\tau_a}
\frac{\partial\phi_m}{\partial\tau_b}, a,b = \overline{1,k},
m=\overline{1,\infty}$,
$$A^I = i(E+M^I)^{-1}(E-M^I)$$
and from the property $M_{mn}\partial\phi^{I*}_m/
\partial\tau_{a} = 0$.

The proof of formula (\ref{e17}) bases on the following Lemma.

\begin{lemma}\label{l3}

Let $y^a,z^a\in l^2,a=\overline{1,k}$, and let $R$
be the operator in $l^2$ of the form $R\kappa = \kappa
-\sum\limits_{c=1}^{k}y^c(z^c,\kappa), \kappa\in l^2$.

Then $\det R = \det(\delta_{ab}-(y^a,z^b)), a,b = \overline{1,k}$.

\end{lemma}

\sbp{Proof.} Choose an orthonormal basis $\{e_s,s=\overline{1,\infty}\}$
in $l^2$ such that only the first $k$ components of the vectors
$y^a, a=\overline{1,k}$ differ from zero. Then
$$\det R = \det(\delta_{ij}-\sum\limits_{c=1}^{\infty}
(e_i,y^c)(z^c,e_j)),\quad i,j=\overline{1,k}.$$

Since the functions $\det (\delta_{ab}-\alpha(y^a,z^b))$
and $\det(\delta_{ij} - \alpha\sum_{c=1}^{\infty}(e_i,y^c)
(z^c,e_j))$ are polynomials in $\alpha$, and
Taylor-series expansions of their logarithms as $\alpha
\to 0$ coinside these functions coinside as well.

This completes the proof of Lemma \ref{l3}.

Lemma \ref{l3} implies that
$$\begin{array}{c}
\det\left(\delta_{mn}+(E+M)_{ms}^{-1}
\frac{\partial\phi_s}{\partial\tau_a}
W_{ab}
\frac{\partial\phi_n}{\partial\tau_b}\right) = \\
= \det\left(\delta_{ab} +
W_{ac}
\frac{\partial\phi_m}{\partial\tau_c}(E+M)_{mn}^{-1}
\frac{\partial\phi_n}{\partial\tau_b}\right),
\quad a,b,c =\overline{1,k}, \quad s,m,n =\overline{1,\infty}
\end{array}$$

Formula (\ref{e17}) now follows. The statement is proved.

Therefore we have
\begin{equation}\label{e18}
\begin{array}{c}
(\K_{[\Lambda^k,r],\tau^{(0)}}^{\epsilon}f_{\alpha})_{I_{\alpha}}
(q_1^{I_{\alpha}},q_2^{I_{\alpha}},\dots) =
\int_{\Lambda^k}\frac{d\tau f(\tau)e^{\pi il_{\alpha}/4}}
{(2\pi)^{k/2}\epsilon^{k/4}} \\
\cdot\frac{\exp\left(\frac{i}{2}\Arg\det\frac{C-iB}{\sqrt{2}}\right)}
{\sqrt{\det(C^{I_{\alpha}}\sqrt{2})}}
\sqrt{\det\frac{1}{i}\left(
\frac{\partial Q_m^{I_{\alpha}}}{\partial\tau_a}(\tau)
A_{mn}^{I_{\alpha}}(\tau)
\frac{\partial Q_n^{I_{\alpha}}}{\partial\tau_b}(\tau) -
\frac{\partial P_m^{I_{\alpha}}}{\partial\tau_a}(\tau)
\frac{\partial Q_m^{I_{\alpha}}}{\partial\tau_b}(\tau)\right)} \\
\cdot\exp\left(\frac{i}{\epsilon}\int\limits_{\tau^{(0)}}^{\tau}
P_m^{I_{\alpha}}(\tau^{'})dQ_m^{I_{\alpha}}(\tau^{'})
+\frac{i}{\epsilon} P_m^{I_{\alpha}}(\tau)
(q_m^{I_{\alpha}}-Q_m^{I_{\alpha}}(\tau))\right. \\
\left. +\frac{i}{2\epsilon}(q_m^{I_{\alpha}}
-Q_m^{I_{\alpha}}(\tau))A_{mn}^{I_{\alpha}}(\tau)
(q_n^{I_{\alpha}}-Q_n^{I_{\alpha}}(\tau))
-\frac{i}{\epsilon}\sum\limits_{j\in {I_{\alpha}}}P_j(\tau^{(0)})
Q_j(\tau^{(0)})\right), \\
m,n =\overline{1,\infty}, a,b = \overline{1,k}
\end{array}
\end{equation}

We choose the sign of $\sqrt{\det(C^{I_{\alpha}}\sqrt{2})}$
so that
$$Re\sqrt{\det(C^{I_{\alpha}}\sqrt{2})}/
\left(e^{\pi il_{\alpha}/4}
\sqrt{\det((C-iB)/\sqrt{2})}\right)>0$$
(this real part cannot vanish since
$$Re\det(C^{I_{\alpha}}\sqrt{2})/(i^{l_{\alpha}}
\det((C-iB)/\sqrt{2}))>0).$$

Find the asymptotics of integral (\ref{e18}). Since $\Im A^I>0$
only the domain in $\R^k$ where
$(q_m^{I_{\alpha}}-Q_m^{I_{\alpha}}(\tau))
(q_m^{I_{\alpha}}-Q_m^{I_{\alpha}}(\tau))\sim\epsilon$
provides a non exponentially
small contribution. Denote by $\tilde\tau(q)$ the point of
minimum of
$(q_m^{I_{\alpha}}-Q_m^{I_{\alpha}}(\tau))
(q_m^{I_{\alpha}}-Q_m^{I_{\alpha}}(\tau))$.

Expanding the exponent in a neighborhood of the point
$\tilde\tau$ we obtain
$$\begin{array}{c}
\int_{\tau^{(0)}}^{\tau}
P_m^{I_{\alpha}}(\tau^{'})
dQ_m^{I_{\alpha}}(\tau^{'}) +
P_m^{I_{\alpha}}(\tau)
(q_m^{I_{\alpha}}-
Q_m^{I_{\alpha}}(\tau)) +
\frac{1}{2}
(q_m^{I_{\alpha}}-
Q_m^{I_{\alpha}}(\tau))A_{mn}^{I_{\alpha}}(\tau)
(q_n^{I_{\alpha}}-
Q_n^{I_{\alpha}}(\tau)) \\
=
\int_{\tau^{(0)}}^{\tilde\tau}
P_m^{I_{\alpha}}(\tau^{'})
dQ_m^{I_{\alpha}}(\tau^{'}) +
P_m^{I_{\alpha}}(\tilde\tau)\xi_m\sqrt{\epsilon}
+\frac{\epsilon}{2}\xi_mA_{mn}^{I_{\alpha}}(\tilde\tau)\xi_n
-\frac{\epsilon}{2}
\frac{\partial P_m^{I_{\alpha}}}{\partial\tilde\tau_a}(\tilde\tau)
\frac{\partial Q_m^{I_{\alpha}}}{\partial\tilde\tau_b}(\tilde\tau)t_at_b \\
+\epsilon\left(
\frac{\partial P_m^{I_{\alpha}}}{\partial\tilde\tau_a}(\tilde\tau) -
\frac{\partial Q_n^{I_{\alpha}}}{\partial\tilde\tau_a}(\tilde\tau)
A_{mn}^{I_{\alpha}}(\tilde\tau)\right)t_a\xi_m
+\frac{\epsilon}{2}
\frac{\partial Q_m^{I_{\alpha}}}{\partial\tilde\tau_a}(\tilde\tau)
A_{mn}^{I_{\alpha}}(\tilde\tau)
\frac{\partial Q_n^{I_{\alpha}}}{\partial\tilde\tau_b}(\tilde\tau)t_at_b
+\dots \\
m,n = \overline{1,\infty},\quad a,b=\overline{1,k}, \quad
t_a=(\tau_a-\tilde\tau_a)/\sqrt{\epsilon},\quad
\xi_m = (q_m-Q_m(\tilde\tau))/\sqrt{\epsilon}
\end{array}$$

Integrating over $t_a,\xi_m$ we obtain the following asymptotics

\begin{equation}\label{e19}
\begin{array}{c}
(\K_{[\Lambda^k,r],\tau^{(0)}}^{\epsilon}f_{\alpha})_{I_{\alpha}}
(q_1^{I_{\alpha}},q_2^{I_{\alpha}},\dots) =
\frac{e^{i\pi l_{\alpha}/4+\frac{i}{2}\Arg\det\left(
\frac{C-iB}{\sqrt{2}}(\tilde\tau)\right)}}
{\sqrt{\det(C^I(\tilde\tau)\sqrt{2})}} \\
\cdot f(\tilde\tau)\epsilon^{k/4}\exp\left\{
\frac{i}{\epsilon}\int\limits_{\tau^{(0)}}^{\tilde\tau}
P_m^I(\tau^{'})dQ_m^I(\tau^{'})
+P_m^I(\tilde\tau)(q_m^I-Q_m^I(\tilde\tau)) \right. \\
\left.+\frac{1}{2}(q_m^I-Q_m^I(\tilde\tau))(B^I(C^I)^{-1})_{mn}
(q_n-Q_n^I(\tilde\tau))
-\frac{i}{\epsilon}\sum\limits_{j\in I}P_j(\tau^{(0)})
Q_j(\tau^{(0)})\right\} \\
m,n =\overline{1,\infty}
\end{array}
\end{equation}
\begin{remark}

It follows therefore that the definition of a canonical operator in
the present paper differs from the classical definition only
by the factor $\const\exp(\frac{i}{2}\Arg\det\frac{C-iB}{\sqrt{2}})$.
This factor changes the quantization conditions and
the transport equation for function $f$.
\end{remark}

This heuristic will not however be used in the proof of the theorem.

\section{Time evolution of a Lagrangian manifold with complex germ}

Consider now a canonical transformation of a Lagrangian manifold with
complex germ.

Let
\begin{equation}\label{e20}
H(\phi^*,\phi) =
\sum\limits_{m=1}^s
\sum\limits_{n=1}^s
H_{i_1\dots i_m j_1\dots j_n}^{(m,n)}
\phi_{i_1}^*\dots\phi_{i_m}^*\phi_{j_1}\dots\phi_{j_m}
\end{equation}
$$ \bar
H_{i_1\dots i_m j_1\dots j_n}^{(m,n)} =
H_{j_1\dots j_n i_1\dots i_m}^{(n,m)} $$
$H_{i_1\dots i_mj_1\dots j_n}^{(m,n)}$ is symmetric separately
over $i_1,\dots,i_m$ and over $j_1,\dots,j_n$.

\begin{definition}

We say that a canonical transformation $\D_H^t,t\in[0,T]$
of a Lagrangian manifold with complex germ $[\Lambda_t^k,r_t]$
{\em corresponds to the Hamiltonian} $H$, if
\begin{enumerate}
\item there exists
 on the segment $[0,T]$
a solution $\phi_j(\tau,t)$
of the Cauchy problem
\begin{eqnarray}\label{e21}
i\dot\phi_j & = & \frac{\partial H}{\partial\phi_j^*}(\phi_j^*,\phi),
\nonumber \\
\phi_j(\tau,t)\vert_{t=0} & = & \phi_j^0(\tau)
=(Q_j^0(\tau)+iP_j^0(\tau))/\sqrt{2}, \\
& & (P^0(\tau),Q^0(\tau))\in\Lambda_0^k \nonumber
\end{eqnarray}
such that for any
$\lambda_1,\dots,\lambda_k,\quad \lambda_i\in
\{0,1,2,\dots\},\quad i=\overline{1,k} $
the derivatives
$$\frac{\partial^{\lambda_1+\dots +\lambda_k}}
{\partial\tau_1^{\lambda_1}\dots\partial\tau_k^{\lambda_k}}
\phi_j(\tau,t)$$
exist, and the series $\sum\limits_{j=1}^{\infty}
\vert\frac{\partial^{\lambda_1+\dots+\lambda_k}}
{\partial\tau_1^{\lambda_1}\dots\partial\tau_k^{\lambda_k}}
\phi_j(\tau,t)\vert^2$ converges.

\item there exists
 on the segment $[0,T]$ a solution of the following Cauchy problem
\begin{eqnarray}\label{e22}
i\dot\Pi_{mn}(\tau,t) & = &
-\frac{\partial^2H}{\partial\phi_m\partial\phi_l^*}
\Pi_{ln}(\tau,t)
-\frac{\partial^2H}{\partial\phi_m\partial\phi_l}
\Omega_{ln}(\tau,t) \nonumber \\
i\dot\Omega_{mn}(\tau,t) & = &
\frac{\partial^2H}{\partial\phi_m^*\partial\phi_l^*}
\Pi_{ln}(\tau,t)
+\frac{\partial^2H}{\partial\phi_m^*\partial\phi_l}
\Omega_{ln}(\tau,t)  \\
& & \Pi_{mn}(\tau,0) =\delta_{mn},\quad \Omega_{mn}(\tau,0) = 0
\quad m,n,l = \overline{1,\infty} \nonumber
\end{eqnarray}
(the arguments $\phi^*(\tau,t),\phi(\tau,t)$ at the
derivatives of $H$ are omitted) such that
for any
$\lambda_1,\dots,\lambda_k,\quad \lambda_i\in
\{0,1,2,\dots\},\quad i=\overline{1,k} $
the derivatives
$$\Pi_{mn}^{(\lambda_1,\dots,\lambda_k)}
= \frac{\partial^{\lambda_1+\dots +\lambda_k}}
{\partial\tau_1^{\lambda_1}\dots\partial\tau_k^{\lambda_k}}
\Pi_{mn}(\tau,t)$$
$$\Omega_{mn}^{(\lambda_1,\dots,\lambda_k)}
= \frac{\partial^{\lambda_1+\dots +\lambda_k}}
{\partial\tau_1^{\lambda_1}\dots\partial\tau_k^{\lambda_k}}
\Omega_{mn}(\tau,t)$$
exist, the operators $\Pi^{(\lambda_1,\dots,\lambda_k)}$ are
bounded, and the operators $\Omega^{(\lambda_1,\dots,\lambda_k)}$
are Hilbert-Schmidt operators.
\end{enumerate}
\end{definition}

{\em A canonical transformation $\D_H^t$ of a Lagrangian manifold
with complex germ} is a set of transformations
mapping the manifold $\Lambda_0^k$ in the manifold
$$g_H^t\Lambda_0^k=\Lambda_t^k=
\left\{
P_j=\frac{\phi_j(\tau,t)-\phi_j^*(\tau,t)}{\sqrt{2}i},
Q_j=\frac{\phi_j(\tau,t)+\phi_j^*(\tau,t)}{\sqrt{2}}\right\}$$
where $\phi_j(\tau,t)$ is a solution of the system (\ref{e21}),
and the matrices $B(\tau,0), C(\tau,0)$ into the
matrices $B(\tau,t),C(\tau,t)$ such that
\begin{equation}\label{e23}
\left(\begin{array}{c}
C-iB \\ C+iB
\end{array}\right)
(\tau,t) =
\left(\begin{array}{cc}
\Pi & \bar\Omega \\ \Omega & \bar\Pi
\end{array}\right)
(\tau,t)
\left(\begin{array}{c}
C-iB \\ C+iB
\end{array}\right)
(\tau,0)
\end{equation}

\begin{remark}
The substitution
$$\phi_j=(Q_j+iP_j)/\sqrt{2},\phi_j^*=(Q_j-iP_j)/\sqrt{2}$$
makes the equations (\ref{e21}) Hamiltonian with the Hamiltonian function
$H((Q_j-iP_j)/\sqrt{2},(Q_j+iP_j)/\sqrt{2})$.

\end{remark}

\begin{remark}

The equations for the matrices $B$ and $C$ are equations
in variations for this Hamiltonian system (see \cite{M2,BD}).

\end{remark}

\begin{lemma}\label{l4}

The pair consisting of the manifold $\Lambda_t^k$ and
the complex germ corresponding to the matrices $B(\tau,t),
C(\tau,t)$ is a Lagrangian manifold with complex germ.

\end{lemma}

\sbp{Proof.}

First of all,check that the following matrix
\begin{equation}\label{e23*}
\left(\begin{array}{cc}
\Pi & \bar\Omega \\
\Omega & \bar\Pi
\end{array}\right)
\end{equation}
is a matrix of a proper canonical transformation,i.e.
\begin{equation}\label{e24}
\Pi^+\Pi-\Omega^+\Omega=E,\quad\Omega^T\Pi=\Pi^T\Omega,
\Pi\Pi^+-\bar\Omega\Omega^T=E,\Omega\Pi^+=\bar\Pi\Omega^T.
\end{equation}

At the initial moment these properties obviously hold.
Equations (\ref{e22}) imply that
$$(\Pi^+\Pi-\Omega^+\Omega)^{\cdot} = 0, \quad
(\Omega^T\Pi-\Pi^T\Omega)^{\cdot}=0 $$
$$
(\Pi\Pi^+-\bar\Omega\Omega^T)^{\cdot}=0,
(\Omega\Pi^+-\bar\Pi\Omega^T)^{\cdot}=0.
$$
Properties (\ref{e24}) hold therefore at any moment of time.

Check now the axioms of a Lagrangian manifold with complex germ.
Properties (\ref{e24}) imply germ axioms r3, r4 and axiom
m3. Axiom m2 follows from the invertibility of
the matrix (\ref{e23*})
Axioms m1,r5 follow immediately from the definition 3.
Axiom r1 follows from linearity of the system (\ref{e22}).Axiom r2 follows
from the equation (\ref{e21}).

Berezin have shown \cite{B} that
for matrix of a proper canonical transfornation (\ref{e23*})
\newline
$||\Pi^{-1}\bar\Omega||<1$ and operator
$\Pi^{-1}$ is bounded.Formula (\ref{e6}) implies that
$ ||(C+iB)(0,\tau)(C-iB)^{-1}(0,\tau)|| \le  1 $,so
$ ||\Pi^{-1}\bar\Omega(C+iB)(0,\tau)(C-iB)^{-1}(0,\tau)|| < 1 $.
Thus ,the following operator
$$ (C-iB)^{-1}(0,\tau)
(E+\Pi^{-1}\bar\Omega(C+iB)(0,\tau)(C-iB)^{-1}(0,\tau))^{-1}
  \Pi^{-1}(t,\tau)=
 (C-iB)^{-1}(t,\tau)  $$
is bounded,so axiom r6 is also satisfied.
Lemma \ref{l4} is proved.

It is also easy to check that the quantization condition
(\ref{e10}) holds at any moment of time whenever it holds at
the initial moment.

\begin{remark} In the case of a germ in a point the
preserving of the germ aximos has the following sense:
the product of two matrices of proper canonical transformations \cite{B}
$$
\left(\begin{array}{cc}
\Pi & \bar\Omega \\ \Omega & \bar\Pi
\end{array}\right)
\qquad \mbox{and} \qquad
\left(\begin{array}{cc}
G & \bar F \\ F & \bar G
\end{array}\right)
$$
is a matrix of a proper canonical transformation .
\end{remark}

\section{Connection between geometrical and canonical
\newline
quantization and some examples}

Suppose $H(\sqrt{\epsilon}\hat\psi^+,
\sqrt{\epsilon}\hat\psi^-)$
is a selfadjoint operator in $\H$ of the form
\begin{equation}\label{e25}
H(\sqrt{\epsilon}\hat\psi^+,
\sqrt{\epsilon}\hat\psi^-)=
\sum\limits_{m,n=1}^s
H_{i_1\dots i_mj_1\dots j_n}^{(m,n)}
\epsilon^{\frac{m+n}{2}}
\hat\psi_{i_1}^+\dots
\hat\psi_{i_m}^+
\hat\psi_{j_1}^-\dots
\hat\psi_{j_n}^-.
\end{equation}

Consider the Cauchy problem for the equation
\begin{equation}\label{e26}
\begin{array}{rcl}
i\frac{\partial\hat\Phi(t)}{\partial t} & = &
\frac{1}{\epsilon} H(\sqrt{\epsilon}\hat\psi^+,
\sqrt{\epsilon}\hat\psi^-)\hat\Phi(t), \quad
\hat\Phi(t)\in\H, \\
\hat\Phi\vert_{t=0} & = &
\hat\K_{[\Lambda_0^k,r_0],\tau^{(0)}}^{\epsilon}f_0
\end{array}
\end{equation}

Introduce a notation
\begin{equation}\label{e27}
\begin{array}{l}
\hat\Psi^{\epsilon}(t) =
\int\frac{d\tau f(\tau,t)}{(2\pi)^{k/2}\epsilon^{k/4}}
\frac{\sqrt{\det
\frac{\partial\phi_l^*}{\partial\tau_a}(\tau,t)
\frac{\partial\phi_l}{\partial\tau_b}(\tau,t)}}
{\sqrt[4]{\det
 [G^+(\tau,t)G(\tau,t)] }} \\
\cdot\exp\left\{
\frac{1}{\epsilon}[g+\int\limits_{(\tau^{(0)},0)}^{(\tau,t)}
(\phi_l(\tau^{'},t^{'})d\phi_l^*(\tau^{'},t^{'})
-iH(\phi^*(\tau^{'},t^{'}),\phi(\tau^{'},t^{'})dt^{'})]\right\}
\hat\Phi_{\phi(\tau,t),M(\tau,t)} \\
l = \overline{1,\infty}, \quad a,b = \overline{1,k}.
\end{array}
\end{equation}
where
$$
g= \phi_l^*(\tau^{(0)},0)\phi_l(\tau^{(0)},0)/2 +
(\phi_l^*(\tau^{(0)},0)\phi_l^*(\tau^{(0)},0)-
\phi_l(\tau^{(0)},0)\phi_l(\tau^{(0)},0))/4,
$$
\begin{equation}\label{e28}
\begin{array}{c}
f(\tau,t) = \exp\left\{-\frac{i}{4}\int\limits_0^t
dt^{'}
\left(
\frac{\partial^2H}{\partial\phi_m\partial\phi_n}
FG^{-1}\right)_{mn} \right.\\
\left. + \left(
\frac{\partial^2H}{\partial\phi_m^*\partial\phi_n^*}
\overline{FG^{-1}}\right)_{mn} \right\}f_0(\tau)
\end{array}
\end{equation}

\begin{theorem}\label{t1}

Suppose $H_{i_1\dots i_mj_1\dots j_n}^{(m,n)}$ are number sets
such that
$$ \| H(\sqrt{\epsilon}\hat\psi^+,
\sqrt{\epsilon}\hat\psi^-)
\hat\Phi_{\phi(\tau,t),M(\tau,t)}\|
\le C_1,\quad\epsilon\in(0,\epsilon_0),C_1>0.$$
where $M(\tau,t)$ is defined by formula (\ref{e7}). Let the
canonical transformation $\D_H^t$ correspond to
the Hamiltonian $H$, and let the series $\sum\limits_{n,m=1}^{\infty}
\vert\frac{\partial^2H}{\partial\phi_m\partial\phi_n}\vert
^2$ converge.

Then the solution of the Cauchy problem (\ref{e26}) may be
presented in the form
$$\hat\Phi(t)=\hat\Psi^{\epsilon}(t)+\hat\delta_{\epsilon}^{(1)}(t) $$
$$(\hat\delta_{\epsilon}^{(1)}(t),
\hat\delta_{\epsilon}^{(1)}(t)) {\longrightarrow \atop \epsilon
\to 0} 0 $$
$$(\hat\Phi(t),\hat\Phi(t)) = O(1)\quad \mbox{as}\quad
\epsilon \to 0. $$
\end{theorem}

A proof of this theorem will be given below as a corollary
of a more general statement. Hence it will be proved that
the canonical operator really gives the asymptotics
of the Cauchy problem solution, i.e. that the geometrical
quantization is compatible with the canonical one.

Consider now some examples.

\begin{example}\label{p1}

The following approximate solution of the equation (\ref{e26})
coinciding up to a constant factor with the vector $\hat\Phi
_{\phi(t),M(t)}$ corresponds to a zero dimensional isotropic
manifold:
$$\begin{array}{l}
\hat\Psi^{\epsilon}(t) =
\exp(g/\epsilon)
\\
\cdot\frac{
f_0\exp\left(-\frac{i}{4}
\int\limits_0^t
\left(\frac{\partial^2H}{\partial\phi_m\partial\phi_n}
M_{mn}(t^{'})
+\frac{\partial^2H}{\partial\phi_m^*\partial\phi_n^*}
M_{mn}^*(t^{'})\right)dt^{'}\right)}
{\sqrt[4]{\det
 [G^+(t)G(t)] }} \\

\cdot\exp\left(
\frac{1}{\epsilon}\int\limits_{0}^{t}
(\phi_l(t^{'})\dot\phi_l^*(t^{'})
-iH)dt^{'}\right)
\hat\Phi_{\phi(t),M(t)} \\
m,n,l = \overline{1,\infty},
\end{array} $$
where we omit arguments $\phi^*(t^{'}),\phi(t^{'})$ at $H$.

The matrix $M$ satisfies the equation
$$i\dot M_{mn} =
\frac{\partial^2 H}{\partial\phi_m^*\partial\phi_n^*}
+\frac{\partial^2 H}{\partial\phi_m^*\partial\phi_s}M_{sn} $$
$$+M_{ms}\frac{\partial^2 H}{\partial\phi_s\partial\phi_n^*}
+M_{ms}\frac{\partial^2 H}{\partial\phi_s\partial\phi_r}M_{rn},$$
$$m,s,n,r =\overline{1,\infty}$$

In the Fock presentation this vector coincides with the
vector (\ref{e9}) up to a constant factor. Note that all
the components of the vector corresponding to a zero
dimensional isotropic manifold differ from zero.
\end{example}

\begin{example}\label{p2} \cite{MS2,MS3}

Let $H$ be of the form
$$H(\phi^*,\phi) =
\sum\limits_{n=1}^s
H_{i_1\dots i_nj_1\dots j_n}^{(n)}
\phi_{i_1}^*\dots\phi_{i_n}^*\phi_{j_1}\dots\phi_{j_n}.$$

The family $g_H^t$ maps a one dimensional isotropic manifold
of the form
$$\phi_j(\tau)=\tilde\phi_j^0\exp(i\tau),\tau\in[0,2\pi)$$
into an isotropic manifold of the form
$$\phi_j(\tau,t)=\tilde\phi_j(t)\exp(i\tau), \quad \tau\in[0,2\pi),$$
where $\phi_j$ is a solution of the equations (\ref{e21}).

Introduce as usual $F=(C+iB)/\sqrt{2},G=(C-iB)/\sqrt{2}$.

The equations in variations for $F$ and $G$
\begin{eqnarray*}
i\dot F_{mn} & = &
\frac{\partial^2 H}{\partial\phi_m^*\partial\phi_r^*} G_{rn}
+ \frac{\partial^2 H}{\partial\phi_m^*\partial\phi_r} F_{rn} \\
i\dot G_{mn} & = &
- \frac{\partial^2 H}{\partial\phi_m\partial\phi_r^*} G_{rn}
- \frac{\partial^2 H}{\partial\phi_m\partial\phi_r} F_{rn},
\quad m,r,n =\overline{1,\infty}
\end{eqnarray*}
have a solution of the form
$$F(\tau,t) =\tilde F(t)e^{i\tau},\quad
G(\tau,t) =\tilde G(t)e^{-i\tau}.$$

The quantities
$$\det G^+G,\quad \det
\frac{\partial\phi_l^*}{\partial\tau_a}
\frac{\partial\phi_l}{\partial\tau_b}
= \tilde\phi_l^*\tilde\phi_l,\quad l =\overline{1,\infty}$$
do not depend on $\tau$. If $f_0$ does not depend on $\tau$,
the function $f(\tau,t)$ does not depend on $\tau$ as well.

In this case vector (\ref{e27}) has therefore the following form
\begin{eqnarray*}
\hat\Psi^{\epsilon}(t) & = &
\exp(g/\epsilon) \\
& & \times\frac{f(t)\sqrt{\tilde\phi_l^*\tilde\phi_l}}
{\sqrt[4]{\det \tilde G^+(t)\tilde G(t)}}
\exp\left(\frac{1}{\epsilon}\int\limits_0^t
(\tilde\phi_l(t^{'})\dot{\tilde\phi_l}^*(t^{'})-iH)dt^{'}\right) \\
& & \times\int\frac{d\tau}{\sqrt{2\pi}\epsilon^{1/4}}
e^{-\frac{i}{\epsilon}\tau\tilde\phi_l^*\tilde\phi_l
+\frac{1}{\epsilon}\tilde\phi_j(t)(\sqrt{\epsilon}
\hat\psi_j^+e^{i\tau}-\tilde\phi_j^*(t))} \\
& & \cdot\exp\left(\frac{1}{2\epsilon}(\hat\psi_j^+\sqrt{\epsilon}e^{i\tau}
-\phi_j^*(t))\tilde M_{jl}(\hat\psi_l^+\sqrt{\epsilon}e^{i\tau}
-\phi_l^*(t))\right)\hat\Phi_0,\quad j,l = \overline{1,\infty}
\end{eqnarray*}
The quantization condition is of the form
$$\tilde\phi_l^*\tilde\phi_l=\epsilon N, \quad N\in\Z.$$

It is easy to see that integrating over $\tau$ makes zero
all the components of $\hat\Psi^{\epsilon}$ but the $N$-th.
This component has the following form (up to a normalizing factor):
\begin{eqnarray}\label{e29}
(\Psi^{\epsilon})_{i_1\dots i_N}^{(N)} & = &
\tilde\phi_{i_1}\dots\tilde\phi_{i_N}
\sum\limits_{k=0}^{[N/2]}\frac{1}{2^kk!} \nonumber \\
& & \cdot\sum\limits_{1\le j_1\ne\dots\ne j_{2k}\le N}
\frac{
\tilde M_{i_{j_1}i_{j_2}}\dots
\tilde M_{i_{j_{2k-1}}i_{j_{2k}}}}
{\tilde\phi_{i_{j_1}}\dots\tilde\phi_{i_{j_{2k}}}} c(t),
\end{eqnarray}

$$ c(t) =\frac{f(t)}
{\sqrt[4]{detG^+(t)G(t)}}
\exp\{\frac{1}{\epsilon}\int\limits_0^{t}
(\tilde\phi_l(t^{'})\dot{\tilde\phi}_l^*(t^{'})-iH)dt^{'}\}$$

$$\tilde M_{lj} = (FG^{-1})_{lj}-\tilde\phi_l\tilde\phi_j
/(\tilde\phi_m^*\tilde\phi_m), \quad l,j,m =\overline{1,\infty}$$

Formula (\ref{e29}) coincides with that obtained in \cite{MS1}
for the stationary
case. It may be also obtained by taking the $N$-th component
of the state vector corresponding to a germ in a point (see
\cite{MS2,MS3}).

The matrix $M$ that vanishes at the initial moment may differ
from zero at other moments of time. It now follows that
an $N$-particle wave function may not be decomposed into
a product of one-particle wave functions
(i.e. in form (\ref{e29}) with $M=0$) as $\epsilon\to 0,
N\to\infty,\epsilon N\to\const$, even in the case it
is decomposed into such a product at the initial moment.
For the case of classical statistical mechanics this statement
has been stated and proved in \cite{MS2}.

\end{example}

\begin{example}\label{p3}

The approach developped here may be also spread formally
to the case of quantum field theory, for scalar quantum electrodynamics,
for example. Formulas, however, look more simple for the
theory of a real scalar field with self-action. This theory
describes an approximation of $\pi$-mesons interaction \cite{BS,Sw}.
The Lagrangian for the theory is of the form

\begin{equation}\label{e30}
\L = \frac{1}{2}(\partial_{\mu}\phi)(\partial^{\mu}\phi)
-\frac{m^2}{2}\phi^2-\frac{g}{4}\phi^4,
\quad \mu = 0,1,\dots,d-1,
\end{equation}
where $d$ is the space-time dimension.

The following Hamiltonian
\begin{equation}\label{e31}
\begin{array}{rcl}
\H(\hat p(\cdot),\hat q(\cdot)) & = &
\int d^{d-1}{\bf x}(\frac{1}{2}\hat p^2({\bf x}) +
\frac{1}{2}({\bf \nabla}\hat q({\bf x}))^2 \\
& & +\frac{m^2}{2}\hat q^2({\bf x}) +
\frac{g}{4} \hat q^4({\bf x})), \\
& & d^{d-1}{\bf x} = dx_1\dots dx_{d-1}
\end{array}
\end{equation}
corresponds to the Lagrangian (\ref{e30}).

After the substitution $\sqrt{q}\hat p = \hat P, \sqrt{g}
\hat q = \hat Q$ the Hamiltonian (\ref{e31}) becomes
$$\begin{array}{rcl}
\H(\hat P(\cdot),\hat Q(\cdot)) & = &
\frac{1}{g}\int d{\bf x}
(\frac{1}{2}\hat P^2({\bf x})+
\frac{1}{2}({\bf\nabla}\hat Q({\bf x}))^2 \\
& & +\frac{m^2}{2}\hat Q^2({\bf x}) +
\frac{1}{4}\hat Q^4({\bf x})), \\
& & [\hat Q({\bf x}),\hat P({\bf y})]
 =
i g\delta({\bf x} - {\bf y})
\end{array}$$

The theory may be interpreted in the terms of particals
with the help of following creation and annihilation operators:
\begin{eqnarray*}
\hat q({\bf x}) & = &
\frac{1}{(2\pi)^{\frac{d-1}{2}}}
\int \frac{d{\bf l}}{\sqrt{2\sqrt{{\bf l}^2+m^2}}}
\left(\hat\psi^+({\bf l})e^{-i{\bf l}{\bf x}}
+ \hat\psi^-({\bf l})e^{i{\bf l}{\bf x}}\right) \\
\hat p({\bf x}) & = &
\frac{i}{(2\pi)^{\frac{d-1}{2}}}
\int {d{\bf l}}\sqrt{\frac{\sqrt{{\bf l}^2+m^2}}{2}}
\left(\hat\psi^+({\bf l})e^{-i{\bf l}{\bf x}}
- \hat\psi^-({\bf l})e^{i{\bf l}{\bf x}}\right)
\end{eqnarray*}

According to the outlined scheme one may assign solutions of
secondary quantized equations of the sort (\ref{e27}) to each
Lagrangian manifold with complex germ.

For constructing such an asymptotics one must solve a
Hamiltonian system and a system of equations in variations.

The Hamiltonian system has the form
\begin{eqnarray*}
\dot Q({\bf x},t) & = & P({\bf x},t) \\
\dot P({\bf x},t) & = & \Delta Q({\bf x},t)
- m^2 Q({\bf x},t) - Q^3({\bf x},t)
\end{eqnarray*}
and the equations in variations over the variables $\delta P,
\delta Q$ are
\begin{eqnarray*}
\dot B({\bf x},t,\tau) & = &
\Delta C({\bf x},t,\tau)
- m^2 C({\bf x},t,\tau) - 3C({\bf x},t,\tau)
Q^2({\bf x},t,\tau) \\
\dot C({\bf x},t,\tau) & = &
B({\bf x},t,\tau)
\end{eqnarray*}

The outlined conception may be applied without a preliminar
regularization if
$$\int d{\bf x}(
Q({\bf x},t,\tau)\sqrt{m^2-\Delta}
Q({\bf x},t,\tau)
+P({\bf x},t,\tau)\frac{1}{\sqrt{m^2-\Delta}}
P({\bf x},t,\tau)) < \infty, $$
and the operator $BC^{-1}-i\sqrt{m^2-\Delta}$ is a Hilbert-
Schmidt operator.

If these conditions are not satisfied, regularization must
be treated more closely.

The conception of geometrical quantization may be applied
as well to the scalar quantum electrodynamics in
the $\alpha$-gauge with the Lagrangian
\begin{eqnarray*}
\L & = &
-\frac{1}{4}
(\partial_{\mu}A_{\nu}-\partial_{\nu}A_{\mu})
(\partial^{\mu}A^{\nu}-\partial^{\nu}A^{\mu})
-\frac{1}{2\alpha}(\partial_{\mu}A^{\mu})^2 \\
& & +(\partial_{\mu}+ieA_{\mu})\Phi^*
(\partial_{\mu}-ieA_{\mu})\Phi
-m^2\Phi^*\Phi-\frac{g}{4}(\Phi^*\Phi)^2.
\end{eqnarray*}
\end{example}
\section{Complex germ creation and annihilation operators}

The theory of complex germs allows to construct not only
asymptotic solutions of type (\ref{e27}) of secondary quantized
equations. Other asymptotic solutions of the equation (\ref{e26})
may be obtained also with the help of so called germ
creation and annihilation operators.

Suppose a basis on a complex germ is chosen so that the
matrix $A$ from the axiom r1  of a complex germ is
diagonal. Consider following creation and annihilation operators
on a Lagrangian manifold with complex germ:
\begin{equation}\label{e32}
\begin{array}{rcl}
\bar\A_{\alpha}(\tau,t) & = &
\bar G_{m\alpha}(\tau,t)\left(\hat\psi_m^+ -
\frac{\phi_m^*(\tau,t)}{\sqrt{\epsilon}}\right)
-\bar F_{m\alpha}(\tau,t)\left(\hat\psi_m^- -
\frac{\phi_m(\tau,t)}{\sqrt{\epsilon}}\right),\\
\A_{\alpha}(\tau,t) & = &
G_{m\alpha}(\tau,t)\left(\hat\psi_m^- -
\frac{\phi_m(\tau,t)}{\sqrt{\epsilon}}\right)
-F_{m\alpha}(\tau,t)\left(\hat\psi_m^+ -
\frac{\phi_m^*(\tau,t)}{\sqrt{\epsilon}}\right),
\end{array}
\end{equation}
$$F = (C+iB)/\sqrt{2}, \quad
G = (C-iB)/\sqrt{2}, \quad  \alpha =\overline{k+1,\infty},
\quad m = \overline{1,\infty} $$

Let $\nu_{\alpha}\in\{0,1,2,\dots\},\alpha = \overline
{k+1,\infty},\sum\limits_{\alpha=k+1}^{\infty}\nu_{\alpha}
<\infty$.
Let $l(\tau^{'},\tau^{''})$ be a closed path on $\Lambda^k$,
such that it is covered by a path on $M^k$ with the beginning
in point $\tau^{'}$ and the end in point $\tau^{''}$.
Suppose that the following quantization condition holds
for any such path
\begin{equation}\label{e33}
\frac{1}{2\pi\epsilon}\oint_{l(\tau^{'},\tau^{''})}
P_md Q_m = \sum\limits_{\alpha=k+1}^{\infty}
\gamma_{\alpha}(\tau^{'},\tau^{''})\nu_{\alpha}+n,
\quad n\in \Z
\end{equation}
where $\gamma_{\alpha}(\tau^{'},\tau^{''})$ are determined by
the condition $A_{\alpha,\alpha}(\tau^{'},\tau^{''})
= e^{i\gamma_{\alpha}(\tau^{'},\tau^{''})}$.

Consider the following element of $\H$:
\begin{equation}\label{e34}
\begin{array}{rcl}
\hat{\tilde\Psi}^{\epsilon}(t) & = &
\int\frac{d\tau f(\tau,t)}{(2\pi)^{k/2}\epsilon^{k/4}}
\frac
{\sqrt{\det
\frac{\partial\phi_l}{\partial{\tau_a}}(\tau,t)
\frac{\partial\phi_l^*}{\partial{\tau_b}}(\tau,t)}}
{\sqrt[4]{\det
G^+(\tau,t)G(\tau,t)}} \\
& & \cdot\exp\left(\frac{i}{\epsilon}
S(\tau,t)\right)
\bar\A_{k+1}^{\nu_{k+1}}
\bar\A_{k+2}^{\nu_{k+2}}\dots
\hat\Phi_{\phi(\tau,t),M(\tau,t)},\\
& & l = \overline{1,\infty},\quad a,b = \overline{1,k}
\end{array}
\end{equation}
where
$$\begin{array}{rcl}
S(\tau,t) & = &
\frac{1}{2i}
\phi_l^*(\tau^{(0)},0)
\phi_l(\tau^{(0)},0)
+\frac{1}{4i}(
\phi_l^*(\tau^{(0)},0)
\phi_l^*(\tau^{(0)},0)
-\phi_l(\tau^{(0)},0) \\
& & \phi_l(\tau^{(0)},0))
+\int_{(\tau^{(0)},0)}^{(\tau,t)}
(\frac{1}{i}
\phi_l(\tau^{'},t^{'})
d\phi_l^*(\tau^{'},t^{'}) \\
& & -H(\phi^*(\tau^{'},t^{'}),\phi(\tau^{'},t^{'}))dt^{'}),
\end{array}$$
and $f(\tau,t), M(\tau,t)$ are determined from the formulas
(\ref{e28}) and (\ref{e7}) respectively. Let $\hat\Phi(t)$ be
the solution of (\ref{e26}), satisfying the initial condition
$\hat\Phi(0)=\hat{\tilde\Psi}^{\epsilon}(0)$.

\begin{theorem}

Suppose
$$\| H(
\sqrt{\epsilon}\psi^+,
\sqrt{\epsilon}\psi^-)
\bar\A_{k+1}^{\nu_{k+1}}
\bar\A_{k+2}^{\nu_{k+2}}\dots\hat\Phi_{\phi(\tau,t),M(\tau,t)}\|
\le C_2,\quad\epsilon\in(0,\epsilon_0),C_{2}>0.$$

Let the canonical transformation $\D_H^t$ correspond to the
Hamiltonian $H$, and the series $\sum\limits_{m,n=1}^{\infty}
\vert\partial^2H/\partial\phi_m\partial\phi_n\vert^2$ converge.
Then
$$\hat\Phi(t)=\hat{\tilde\Psi}^{\epsilon}(t)+\hat\delta_{\epsilon}(t),$$
$$(\hat\delta_{\epsilon}(t),
\hat\delta_{\epsilon}(t)) {\longrightarrow \atop \epsilon\to 0} 0,
\quad (\hat\Phi(t),\hat\Phi(t))=O(1), \quad \epsilon\to 0.$$
\end{theorem}

\begin{remark}
In the case of C-Lagrangian manifold  \cite{M1}
,when $\gamma_{\alpha}=0$, a canonical
operator can be generalized.

Let $ f^{(n)}_{\alpha_1,...,\alpha_n}(\tau)$, where
$\alpha_i=\overline{k+1,\infty}, n=\overline{0,\infty} $ be a set of
infinitely differentiable functions
with compact support on the Lagrangian manifold
with complex germ $[\Lambda^k,r]$,where the following series
$$ \sum\limits_{n=1}^{\infty} \sum\limits_{\alpha_1,...,\alpha_n=k+1}
^{\infty} | f^{(n)}_{\alpha_1,...,\alpha_n}(\tau) | ^2 $$
converges at all $\tau$.

The following element of $\H$ is assigned to the Lagrangian manifold
with complex germ $[\Lambda^k,r]$, the set of functions
$ f^{(n)}_{\alpha_1,...,\alpha_n}(\tau)$, the point
 $\tau^{(0)}\in\Lambda^k$,
and the number $\epsilon>0$:

\begin{equation}\label{e34*}
\begin{array}{rcl}
\hat\K_{[\Lambda^k,r],\tau^{(0)}}^{\epsilon} f & = &
\int\frac{d\tau }{(2\pi)^{k/2}\epsilon^{k/4}}
\frac
{\sqrt{\det
\frac{\partial\phi_l}{\partial{\tau_a}}(\tau)
\frac{\partial\phi_l^*}{\partial{\tau_b}}(\tau)}}
{\sqrt[4]{\det
G^+(\tau)G(\tau)}} \\
& & \cdot\exp\left(\frac{i}{\epsilon}
S(\tau)\right)
 \sum\limits_{n=1}^{\infty} \sum\limits_{\alpha_1,...,\alpha_n=k+1}
^{\infty}  \frac{1}{\sqrt{n!}}
f^{(n)}_{\alpha_1,...,\alpha_n}(\tau)

\bar\A_{\alpha_1}
\bar\A_{\alpha_2}\dots
\bar\A_{\alpha_n}

\hat\Phi_{\phi(\tau),M(\tau)},\\
& & l = \overline{1,\infty},\quad a,b = \overline{1,k}
\end{array}
\end{equation}

If $f^{(0)}=f$ and other $f^{(n)}=0$ then canonical operator
(\ref{e34*}) is equal to canonical operator (\ref{e11}).

\end{remark}

\begin{example}\label{p41}

Consider the Hamiltonian from example \ref{p2}, and suppose
$\phi_l(t)$ is a solution of the form $\phi_l(t)=\tilde\phi_l
\exp(-i\Omega t)$ of equation (\ref{e21}).

Suppose that there exist matrices $\tilde F$ and $\tilde G$
such that
\begin{eqnarray*}
(\beta_l+\Omega)\tilde G_{ml} & = &
\sum\limits_{r=1}^{\infty}\frac{\partial^2 H}
{\partial\phi_m\partial\phi_r^*}\tilde G_{rl}
+\sum\limits_{r=1}^{\infty}\frac{\partial^2 H}
{\partial\phi_m\partial\phi_r}\tilde F_{rl}, \\
-(\beta_l-\Omega)\tilde F_{ml} & = &
\sum\limits_{r=1}^{\infty}\frac{\partial^2 H}
{\partial\phi_m^*\partial\phi_r^*}\tilde G_{rl}
+\sum\limits_{r=1}^{\infty}\frac{\partial^2 H}
{\partial\phi_m^*\partial\phi_r}\tilde F_{rl},
\end{eqnarray*}
$$\beta_l\in\R,\qquad \tilde G_{m1}=\tilde\phi_m^*,
\qquad\tilde F_{m1}=-\tilde\phi_m$$
$$\tilde F^T\tilde G = \tilde G^T\tilde F,\quad
\tilde G_{m\alpha}^*\tilde G_{m\beta}
-\tilde F_{m\alpha}^*\tilde F_{m\beta}
=\delta_{\alpha\beta},\qquad m,l=\overline{1,\infty},
\alpha,\quad \beta=\overline{2,\infty}$$
 $\tilde F$ is a Hilbert--Schmidt operator
and $\tilde G$ has a bounded inverse operator.

Consider a Lagrangian manifold with complex germ corresponding to
the following function $\phi_l$
$$\phi_l(\tau,t)=\tilde\phi_l\exp(i(\tau-\Omega t))$$
with matrices $F$ and $G$
\begin{eqnarray*}
F_{ml}(\tau,t) & = &
\tilde F_{ml}\exp(-i(\beta_l/\Omega-1)(\tau-\Omega t)), \\
G_{ml}(\tau,t) & = &
\tilde G_{ml}\exp(-i(\beta_l/\Omega+1)(\tau-\Omega t)).
\end{eqnarray*}
The matrix $A(0,2\pi)$ is
$$A(0,2\pi) = \diag\{1,\exp(2\pi i\beta_2/\Omega),
\exp(2\pi i\beta_3/\Omega),\dots\}.$$
And the quantization condition is
$$\tilde\phi_l^*\tilde\phi_l=\epsilon N+\epsilon\sum\limits
_{\lambda=2}^{\infty}\beta_{\lambda}\nu_{\lambda}/\Omega.$$
Vector (\ref{e34}) depends on time as $\exp(-i\E t)$, where
$$\E=\frac{1}{\epsilon}H(\phi^*,\phi)+\frac{1}{4}
\sum\limits_{m,n=1}^{\infty}
\left(
\frac{\partial^2H}{\partial\phi_m\partial\phi_n}
(FG^{-1})_{mn}
+\frac{\partial^2H}{\partial\phi_m^*\partial\phi_n^*}
(FG^{-1})_{mn}^*\right).$$
Hence it is an approximate solution of the stationary equation
\begin{equation}\label{e35}
\E\hat\Phi = \frac{1}{\epsilon}H(
\hat\psi^+\sqrt{\epsilon},
\hat\psi^-\sqrt{\epsilon})\hat\Phi
\end{equation}

\end{example}

\begin{example}\label{p42}
Consider now the complex germ given by the following matrices:
\begin{eqnarray*}
F_{ml}(\tau,t) & = & F_{ml}\exp(-i\Omega t+i\tau)
\exp(i\beta_lt) \\
G_{ml}(\tau,t) & = & G_{ml}\exp(i\Omega t-i\tau)
\exp(i\beta_lt)
\end{eqnarray*}
The matrix $A(0,2\pi)$ is the identity matrix, and the quantization
condition has the form
$$\tilde\phi_l^*\tilde\phi_l=\epsilon N.$$
In this example vector (\ref{e34}) is also an approximate solution of the
equation (\ref{e35}) for
$$\E = \frac{1}{\epsilon}H(\phi^*,\phi)+
\frac{1}{4}\left(
\frac{\partial^2H}{\partial\phi_m\partial\phi_n}(FG^{-1})_{mn}
+\frac{\partial^2H}{\partial\phi_m^*\partial\phi_n^*}
(FG^{-1})_{mn}^*\right)+\sum\limits_{\lambda=2}^{\infty}
\beta_\lambda\nu_\lambda,\qquad
m,n=\overline{1,\infty}$$
This solution coinside with that given in \cite{MS1}.

Find now a relation between the results of examples 4 and 5.
 Let $\phi_{(1)_m}, \phi_{(2)_m}$ be solutions of the equations
$$\Omega_1\phi_{(1)_m} = \partial H/\partial\phi_{(1)_m}^*,\quad
\Omega_2\phi_{(2)_m} = \partial H/\partial\phi_{(2)_m}^*$$
such that
$$\phi_{(1)_m}^*\phi_{(1)_m}=\epsilon N+\epsilon\sum\limits_{\lambda=1}
^{\infty}\beta_{\lambda}\nu_{\lambda}/\Omega_1,\qquad
\phi_{(2)_m}^*\phi_{(2)_m}=\epsilon N, $$
$$\phi_{(1)_m}=\phi_{(2)_m}+\chi_m,\qquad \|\chi\|=O(\epsilon),\qquad
\Omega_2-\Omega_1=O(\epsilon).$$
We have
$$\begin{array}{rcl}
H(\phi_{(1)}^*,\phi_{(1)}) -
H(\phi_{(2)}^*,\phi_{(2)}) & = &
\Omega_2(\phi_{(2)_m}\chi_m^*+\phi_{(2)_m}^*\chi_m)+O(\epsilon^2),\\
(\phi_{(1)},\phi_{(1)}) -
(\phi_{(2)},\phi_{(2)}) & = &
\phi_{(2)_m}\chi_m^*+\phi_{(2)_m}^*\chi_m+O(\epsilon^2),
\quad m=\overline{1,\infty}
\end{array}$$
It now follows that the values of $\E$ in two examples
coincide up to $O(\epsilon^1)$.

\end{example}

\begin{remark}

We have shown for a specific example how different quantization
conditions may be used with simultanious change of the
transport equation. Conviniency determines the choice of
quantization conditions. In the finite dimensional case
a quantization condition of the kind
$$\frac{1}{2\pi\epsilon}\oint_{l(\tau^{'},\tau^{''})}
P_mdQ_m=\sum\limits_{\alpha=k+1}^{D}
\gamma_{\alpha}(\tau^{'},\tau^{''})
(\nu_{\alpha}+\frac{1}{2})+n,\quad n\in\Z$$
is often used (\cite{BD}).

In the infinite dimensional case this condition makes no
sense even for the simplest Hamiltonian $H(\phi^*,\phi) =
\sum\limits_{i=1}^{\infty}\phi_i^*\phi_i$ and for the
isotropic manifold with complex germ as in example 4.1.

\end{remark}

\begin{example}\label{p43}

Consider one more (heuristic) way to deduce a formula for $\E$
(see \cite{K}).

Note, that if a family $\D_H^t$ acts on a two dimensional
isotropic manifold with complex germ in the following way
$$\begin{array}{rcl}
(g^t\phi)(\tau_1,\tau_2)& =  &
\phi(\tau_1+\Omega_1t,\tau_2+\Omega_2t) \\
F(\tau_1,\tau_2,t) & = &
F(\tau_1+\Omega_1t,\tau_2+\Omega_2t) \\
G(\tau_1,\tau_2,t) & = &
G(\tau_1+\Omega_1t,\tau_2+\Omega_2t)
\end{array}$$
then multiplication of the function $f\tau_1,\tau_2,0)$
by $\exp(i(n_1\tau_1+n_2\tau_2))$ leads to multiplication
of the vector (\ref{e34}) by $\exp(-i(n_1\Omega_1+n_2\Omega_2)t),
\quad n_1,n_2\in\Z$.

Consider now "almost invariant" two dimensional isotropic manifolds
close to one dimensional manifolds in examples 4,5:
$$\phi_l(t,\tau_1,\tau_2)=e^{-i(\tau_1-\Omega t)}
(\phi_l+\delta(
\tilde F_{lm}e^{i(\tau_2+\beta^mt)}
+ \tilde G_{lm}^*e^{-i(\tau_2+\beta^mt)})),
\qquad\delta\to 0. $$
The values of $\E$ corresponding to these almost invariant
manifolds are
\newline
$\E_m^{(\mu)}=\E_m^{(0)}+\beta_m\mu_m,\quad
\mu_m\in\Z$. Since these values must agree with the value for
a neighborhood of a one dimensional manifold the
last value is of the form $\E^{(0)}+\sum_{m=1}^{\infty}
\beta_m\nu_m,\quad \nu_m\in\{0,1,2,\dots\}$.

\end{example}
\section{Proof of the theorem}

Let us now prove theorem 2 on the base of the following lemma.

Introduce notation
$$\begin{array}{rcl}
S(\tau) & = &
\frac{1}{2i}
\phi_j^*(\tau^{(0)})\phi_j(\tau^{(0)})
+\frac{1}{4i}(
\phi_j^*(\tau^{(0)})\phi_j^*(\tau^{(0)})
-\phi_j(\tau^{(0)})\phi_j(\tau^{(0)})) \\
& & +\frac{1}{i}\int_{\tau^{(0)}}^{\tau}
\phi_j(\tau^{'})d\phi_j^*(\tau^{'}),\quad j=\overline{1,\infty}
\end{array}$$

We assign to number sets $\D_{ij}(\tau),\quad i,j=\overline{1,\infty}$
and $(Y_a^n)_{i_1\dots i_n}(\tau),\quad a=1,2, \quad n=\overline{1,s},
\quad i_p=\overline{1,\infty}$ the following operators in $\H$:
$$\begin{array}{rcl}
Y_a\left(\tau,\hat\psi^+-\frac{\phi^*(\tau)}{\sqrt{\epsilon}}\right) & = &
\sum\limits_{n=0}^s(Y_a^n)_{i_1\dots i_n}(\tau)
(\hat\psi_{i_1}^+-\frac{\phi_{i_1}^*(\tau)}{\sqrt{\epsilon}})
\dots(\hat\psi_{i_n}^+-\frac{\phi_{i_n}^*(\tau)}{\sqrt{\epsilon}}) \\
& & \times\exp\{\frac{1}{2}
(\hat\psi_{i}^+-\frac{\phi_{i}^*(\tau)}{\sqrt{\epsilon}})
\D_{ij}(\tau)
(\hat\psi_{i}^+-\frac{\phi_{i}^*(\tau)}{\sqrt{\epsilon}})\}, \\
& & i,j,i_1,\dots,i_n=\overline{1,\infty},\qquad a=1,2.
\end{array} $$

Introduce also the following notation
$$\begin{array}{rcl}
\hat\Phi_a^{\epsilon} & = &
\int\limits_{\Lambda^k}
\frac{d\tau e^{\frac{i}{\epsilon}S(\tau)}}{\epsilon^{k/4}}
Y_a(\tau,
\hat\psi^+-\frac{\phi^*(\tau)}{\sqrt{\epsilon}}) \\
& & \cdot \exp(\frac{1}{\epsilon}
\phi_j(\tau)
(\hat\psi_{j}^+\sqrt{\epsilon}-\phi_{j}^*(\tau)))
\hat\Phi_0,\qquad a=1,2,\quad j=\overline{1,\infty}.
\end{array}$$
The commutation relations for the operators $\psi_j^{\pm}$ imply
$$\begin{array}{rcl}
\hat\Phi_a^{\epsilon} & = &
\int\frac{d\tau\exp(\frac{i}{\epsilon}S(\tau)
-\frac{1}{2\epsilon}\phi_j(\tau)\phi_j^*(\tau))}
{\epsilon^{k/4}} \\
& & \times e^{\frac{1}{\sqrt{\epsilon}}
(\phi_j(\tau)\hat\psi_j^+-\phi_j^*(\tau)\hat\psi_j^-)}
Y_a(\tau,\hat\psi^+)\hat\Phi_0,\qquad j=\overline{1,\infty},\quad a=1,2
\end{array}$$

\begin{lemma}\label{l5}
Let the matrix $\D_{ij}(\tau)$ correspond to a Hilbert--Schmidt
operator $\D$ in $l^2$, $\|\D\|<1$, $\D_{ij}(\tau)$ and
$(Y_a^n)_{i_1\dots i_n}(\tau)$ being smooth functions of $\tau$, functions
$(Y_a^n)_{i_1\dots i_n}(\tau)$
having compact supports,
$$\sum\limits_{i_1\dots i_n}^{\infty}
\vert(Y_a^n)_{i_1\dots i_n}(\tau)\vert^2 < c,\qquad n =\overline{1,s},
\quad a=1,2. $$
Then
\begin{eqnarray}\label{e36}
(\hat\Phi_1^{\epsilon},
\hat\Phi_2^{\epsilon}) & {\longrightarrow \atop \epsilon\to 0} &
\int_{\Lambda^k}d\tau\int_{\R^k}d\xi(
Y_1(\tau,\hat\psi^+)\hat\Phi_0, \nonumber \\
& & \exp(\xi_b(
\frac{\partial\phi_j}{\partial\tau_b}(\tau)\hat\psi_j^+
-\frac{\partial\phi_j^*}{\partial\tau_b}(\tau)\hat\psi_j^-))
Y_2(\tau,\hat\psi^+)\hat\Phi_0),  \\
& & b =\overline{1,k},\qquad \xi\in\R^k,\qquad j=\overline{1,\infty}.
\nonumber \end{eqnarray}
\end{lemma}

\sbp{Proof.}

We have
$$\begin{array}{rcl}
(\hat\Phi_1^{\epsilon},
\hat\Phi_2^{\epsilon}) & = &
\int\frac{d\tau d\tau^{'}}{\epsilon^{k/2}}
e^{\frac{i}{\epsilon}(S(\tau^{'})-S^*(\tau))
-\frac{1}{2\epsilon}
(\phi_j(\tau)\phi_j^*(\tau)
+ \phi_j(\tau^{'})\phi_j^*(\tau^{'})}) \\
& & \times e^{\frac{1}{2\epsilon}
(\phi_j^*(\tau)\phi_j(\tau^{'}
- \phi_j^*(\tau^{'})\phi_j(\tau))}
(\hat\Phi_0,Y_1^*(\tau,\hat\psi^-)\\
& & \cdot e^{\frac{1}{\sqrt{\epsilon}}
(\hat\psi_j^+(\phi_j(\tau^{'})-\phi_j(\tau))
-\hat\psi_j^-(\phi_j^*(\tau^{'})-\phi_j^*(\tau))}
Y_2(\tau^{'},\hat\psi^+)\hat\Phi_0),\qquad
j=\overline{1,\infty}\end{array}$$

It is easy to show that the contribution of the integrating domain
$(\tau_b-\tau_b^{'})(\tau_b-\tau_b^{'})>\delta_1 \epsilon^{1/2-\lambda} ,
\quad\delta_1>0,\lambda>0,\quad b
= \overline{1,k}$ in the integral (\ref{e36}) is exponentially small.
For
$(\tau_b-\tau_b^{'})(\tau_b-\tau_b^{'})\le\delta_1 \epsilon^{1/2-\lambda}$
the integrand may be presented as $$\begin{array}{l}
\frac{1}{\epsilon^{k/2}}e^{\sqrt{\epsilon}\xi_b\xi_c\xi_d
R_{bcd}(\tau,\tau^{'})}(\hat\Phi_0,Y_1^*(\tau,\hat\psi^-) \\
\times e^{\xi_b(
\hat\psi_j^+\frac{\partial\phi_j}{\partial\tau_b}
-\hat\psi_j^-\frac{\partial\phi_j^*}{\partial\tau_b})
+\sqrt{\epsilon}(\hat\psi_j^+J_j(\tau,\tau^{'})-
\hat\psi_j^-J_j(\tau,\tau^{'}))} \\
\times Y_2(\tau^{'},\hat\psi^+)\hat\Phi_0), \quad j=\overline{1,\infty},
\quad b,c,d=\overline{1,k}, \\
\xi_b=(\tau_b^{'}-\tau_b)/\sqrt{\epsilon},\quad
\vert R_{bcd}(\tau,\tau^{'})\vert<\const, \\
\sum\limits_{j=1}^{\infty}\vert J_j(\tau,\tau^{'})\vert^2<\const.
\end{array}$$

Integrating over $\tau$ and $\xi$ as $\epsilon \to 0$ we obtain
(\ref{e36}). Lemma \ref{l5}  is proved.

Lemma \ref{l5} implies

\begin{corollary}

$\|\hat{\tilde\psi}^{\epsilon}(t)\|=O(1)$ as $\epsilon\to 0$.

\end{corollary}

\sbp{Proof.} It is sufficient to prove that only the integral
\begin{equation}\label{e37}\begin{array}{l}
\int d^k\xi(\hat\Phi_0,
e^{\frac{1}{2}\hat\psi_m^-M_{mn}^*
\hat\psi_n^-}
\hat a_{k+1}^{\nu_{k+1}}
\hat a_{k+2}^{\nu_{k+2}}\dots \\
\times\exp\{\xi_b(
\frac{\partial\phi_m}{\partial\tau_b}\hat\psi_m^+
-\frac{\partial\phi_m^*}{\partial\tau_b}\hat\psi_m^-)\}
\hat a_{k+1}^{+\nu_{k+1}}
\hat a_{k+2}^{+\nu_{k+2}}\dots
e^{\frac{1}{2}\hat\psi_m^+M_{mn}
\hat\psi_n^+}\hat\Phi_0) \\
\hat a_{\alpha}^+=\bar G_{m\alpha}\hat\psi_m^+-
\bar F_{m\alpha}\hat\psi_m^-,\quad
\hat a_{\alpha}= G_{m\alpha}\hat\psi_m^--
 F_{m\alpha}\hat\psi_m^+,\qquad
m,n=\overline{1,\infty},\quad b=\overline{1,k}
\end{array}\end{equation}
differs from zero. Show first that
\begin{equation}\label{e38}\begin{array}{c}
\int d^k\xi(\hat X,\exp(\xi_b
(\frac{\partial\phi_m}{\partial\tau_b}\hat\psi_m^+
-\frac{\partial\phi_m^*}{\partial\tau_b}\hat\psi_m^-))
\hat a_{\beta}\exp(
\frac{1}{2}\hat\psi_m^+M_{mn}\hat\psi_n^+)
\hat\Phi_0)=0,\\
\beta\in\{k+1,k+2,\dots\},\qquad b=\overline{1,k},
\quad m,n=\overline{1,\infty}, \\
\hat X =
\hat a_{\alpha_1}^+\dots
\hat a_{\alpha_r}^+
\hat a_{\beta_1}\dots
\hat a_{\beta_s}
\exp(
\frac{1}{2}\hat\psi_m^+M_{mn}\hat\psi_n^+)\hat\Phi_0, \\
\alpha_1,\dots,\alpha_r,\beta_1,\dots,\beta_s
\in\{k+1,k+2,\dots\}
\end{array}\end{equation}

Indeed, make use of the presentation of the scalar
product in the left hand side
of (\ref{e38})  as a functional integral (see \cite{B}):
$$\int d^k\xi\int\prod dz^* dz
X(z^*)(G_{m\beta}\partial/\partial z_m
-F_{m\beta} z_m) $$
$$\times\exp\{\xi_b
(\frac{\partial\phi_m}{\partial\tau_b}z_m
-\frac{\partial\phi_m}{\partial\tau_b}\frac{\partial}{\partial z_m})\}
\exp(\frac{1}{2}z_mM_{mn}z_n)e^{-z_m^*z_m}$$

Integrating over $\xi$ and using the property $M_{mn}\frac
{\partial\phi_n^*}{\partial\tau_b}=0$ we conclude that the
functional integral is equal to zero.

Check now that integral (\ref{e37}) differs from zero. The commutation
relations
$$[\hat a_{\alpha},\hat a_{\beta}]=
[\hat a_{\alpha}^+,\hat a_{\beta}^+]=
[\hat a_{\alpha},
\frac{\partial\phi_m}{\partial\tau_b}\hat\psi_m^+
- \frac{\partial\phi_m^*}{\partial\tau_b}\hat\psi_m^-] $$
$$[\hat a_{\alpha}^+,
\frac{\partial\phi_m}{\partial\tau_b}\hat\psi_m^+
- \frac{\partial\phi_m^*}{\partial\tau_b}\hat\psi_m^-] = 0,\qquad
[\hat a_{\alpha},\hat a_{\beta}^+]=\delta_{\alpha,\beta}$$
make this integral equal to
$$(\nu_{k+1})!
(\nu_{k+2})!\dots
\int d^k\xi(\hat\Phi_0,
\exp(\frac{1}{2}\hat\psi_m^-M_{mn}^*\hat\psi_n^-)$$
$$\times\exp(\xi_b
\frac{\partial\phi_m}{\partial\tau_b}\hat\psi_m^+
- \frac{\partial\phi_m^*}{\partial\tau_b}\hat\psi_m^-))
\exp(\frac{1}{2}\hat\psi_m^+M_{mn}\hat\psi_n^+)\hat\Phi_0.$$

It is easy to check that the last expression differs from zero.
The corollary is proved.

Now we may give a proof of theorem 2. Check that
\begin{equation}\label{e39}
(i\frac{\partial}{\partial t}-\frac{1}{\epsilon}
H(\sqrt{\epsilon}\hat\psi^+,
\sqrt{\epsilon}\hat\psi^-)\hat{\tilde\Psi}^{\epsilon}(t)
{\longrightarrow \atop \epsilon\to 0} 0
\end{equation}

The vector $\hat{\tilde\Psi}^{\epsilon}(t)$ may be presented as
$$\hat{\tilde\Psi}^{\epsilon}(t) =
\int\frac{d\tau\chi(\tau,t)}{\epsilon^{k/4}}
\exp(\frac{i}{\epsilon}S(\tau,t)-
\frac{1}{\epsilon}\phi_j^*(\tau,t)\phi_j(\tau,t))$$
$$\times
e^{\frac{1}{\sqrt{\epsilon}}\phi_j(\tau,t)\hat\psi_j^+}
e^{-\frac{1}{\sqrt{\epsilon}}\phi_j^*(\tau,t)\hat\psi_j^-}
(\hat a_{k+1}^+)^{\nu_{k+1}}
(\hat a_{k+2}^+)^{\nu_{k+2}}\dots
\hat\Phi_{0,M(\tau,t)},\qquad j=\overline{1,\infty}$$
\begin{equation}\label{e40}
\chi(\tau,t)=\frac{f(\tau,t)}{(2\pi)^{k/2}}
\frac{\sqrt{\det
\frac{\partial\phi_j}{\partial\tau_a}(\tau,t)
\frac{\partial\phi_j^*}{\partial\tau_b}(\tau,t)}}
{\sqrt[4]{\det G^+(\tau,t)G(\tau,t)}}
\end{equation}
$$a,b=\overline{1,k},\quad j=\overline{1,\infty}. $$

Since
$$\begin{array}{rcl}
i\frac{\partial}{\partial t}\left(
e^{-\frac{1}{\epsilon}\phi_j^*\phi_j}
e^{\frac{1}{\sqrt{\epsilon}}\phi_j\hat\psi_j^+}
e^{-\frac{1}{\sqrt{\epsilon}}\phi_j^*\hat\psi_j^-}
\right) & = &
e^{-\frac{1}{\epsilon}\phi_j^*\phi_j}
e^{\frac{1}{\sqrt{\epsilon}}\phi_j^*\hat\psi_j^+}
e^{-\frac{1}{\sqrt{\epsilon}}\phi_j\hat\psi_j^-}\\
& & \times(i\frac{\partial}{\partial t}-\frac{i}{\epsilon}
\dot\phi_j^*\phi_j+\frac{i}{\sqrt{\epsilon}}
(\dot\phi_j\hat\psi_j^+-\dot\phi_j^*\hat\psi_j^-))
\end{array}$$
where $j=\overline{1,\infty}$, arguments $t,\tau$ at $\phi_j^*,
\phi_j$ are omitted, we obtain
$$(i\frac{\partial}{\partial t} - \frac{1}{\epsilon}
H(\sqrt{\epsilon}\hat\psi^+,\sqrt{\epsilon}\hat\psi^-))
\hat{\tilde\Psi}^{\epsilon}(t) =
\int\frac{d\tau}{\epsilon^{k/4}}$$
$$\times\exp(\frac{i}{\epsilon}S(\tau,t)
- \frac{1}{2\epsilon}\phi_j^*(\tau,t)\phi_j(\tau,t)
+ \frac{1}{\sqrt{\epsilon}}(
\phi_j(\tau,t)\hat\psi_j^+
-\phi_j^*(\tau,t)\hat\psi_j^-)) $$
$$\times(i\frac{\partial}{\partial t}+\frac{i}{\sqrt{\epsilon}}
(\dot\phi_j(\tau,t)\hat\psi_j^+
-\dot\phi_j^*(\tau,t)\hat\psi_j^-)-\frac{1}{\epsilon}
\frac{\partial S}{\partial t} $$
$$-\frac{i}{\epsilon}
\dot\phi_j^*(\tau,t)\phi_j(\tau,t)
-\frac{1}{\epsilon} H(\phi^*+\sqrt{\epsilon}\hat\psi^+,
\phi+\sqrt{\epsilon}\hat\psi^-))$$
\begin{equation}\label{e41}
\times\chi(\tau,t)
(\hat a_{k+1}^+)^{\nu_{k+1}}
(\hat a_{k+2}^+)^{\nu_{k+2}}\dots
\hat\Phi_{0,M(\tau,t)},\quad j=\overline{1,\infty}.
\end{equation}

The conditions of the theorem imply that the norm of the vector
$$H(
\phi^*(\tau,t)+\sqrt{\epsilon}\hat\psi^+,
\phi(\tau,t)+\sqrt{\epsilon}\hat\psi^-)
(\hat a_{k+1}^+)^{\nu_{k+1}}
(\hat a_{k+2}^+)^{\nu_{k+2}}\dots
\hat\Phi_{0,M(\tau,t)}$$
is uniformly bounded. Since the operator $H(\phi^*+\sqrt{\epsilon}
\hat\psi^+,\phi+\sqrt{\epsilon}\hat\psi^-)$ is a polynomial in
$\sqrt{\epsilon}$ of the form
$$\begin{array}{rcl}
H(\phi^*+\sqrt{\epsilon}\hat\psi^+,
\phi+\sqrt{\epsilon}\hat\psi^-) & = &
\sum\limits_{n=0}^s
\sum\limits_{m=0}^s\epsilon^{\frac{m+n}{2}}
\frac{\partial^{m+n}H}{
\partial\phi_{i_1}^*\dots
\partial\phi_{i_m}^*
\partial\phi_{j_1}\dots
\partial\phi_{j_n}} \\
& & \hat\psi_{i_1}^+
\dots \hat\psi_{i_m}^+
\hat\psi_{j_1}^-
\dots \hat\psi_{j_n}^-,
i_1,\dots,i_m,j_1,\dots,j_n =\overline{1,\infty},
\end{array}$$
the norms of the vectors
$$\frac{\partial^{m+n}H}{
\partial\phi_{i_1}^*\dots
\partial\phi_{i_m}^*
\partial\phi_{j_1}\dots
\partial\phi_{j_n}}
\hat\psi_{i_1}^+
\dots \hat\psi_{i_m}^+
\hat\psi_{j_1}^-
\dots \hat\psi_{j_n}^-
(\hat a_{k+1}^+)^{\nu_{k+1}}
(\hat a_{k+2}^+)^{\nu_{k+2}}\dots
\hat\Phi_{0,M(\tau,t)}$$
are also uniformly bounded. By Lemma \ref{l5} the right hand side
of the formula (\ref{e41}) may be presented as
$$\int\frac{d\tau}{\epsilon^{k/4}}
\exp(\frac{i}{\epsilon}S(\tau,t)
- \frac{1}{2\epsilon}\phi_j^*(\tau,t)\phi_j(\tau,t)$$
$$+ \frac{1}{\sqrt{\epsilon}}(
\phi_j(\tau,t)\hat\psi_j^+
-\phi_j^*(\tau,t)\hat\psi_j^-))
(i\frac{\partial}{\partial t}
-\frac{1}{2}\frac{\partial^2 H}{\partial\phi_m^*\partial\phi_n^*}
\hat\psi_m^+\hat\psi_n^+$$
$$-\frac{1}{2}\frac{\partial^2 H}{\partial\phi_m\partial\phi_n}
\hat\psi_m^-\hat\psi_n^-
-\frac{\partial^2 H}{\partial\phi_m^*\partial\phi_n}
\hat\psi_m^+\hat\psi_n^-
\chi(\tau,t)$$
$$\cdot(\hat a_{k+1}^+)^{\nu_{k+1}}
(\hat a_{k+2}^+)^{\nu_{k+2}}\dots
\hat\Phi_{0,M(\tau,t)} + O(\epsilon^{1/2})$$
$$ j,m,n=\overline{1,\infty}.$$
where we made use of the relations
$$\partial S/\partial t = -i\phi_l\dot\phi_l^*-H(\phi^*,\phi)$$
by the definition of $S$, and
$$i\dot\phi_l = \frac{\partial H}{\partial\phi_l^*},\qquad
-i\dot\phi_l^* = \frac{\partial H}{\partial\phi_l},
\quad l=\overline{1,\infty}$$
by equations (\ref{e21}).

It is easy to check that
\begin{equation}\label{e42}
[i\frac{\partial}{\partial t} - \hat H_2,\hat a_{\alpha}^+] =
[i\frac{\partial}{\partial t} - \hat H_2,\hat a_{\alpha}] =
[i\frac{\partial}{\partial t} - \hat H_2,
\frac{\partial\phi_j}{\partial\tau_a}\hat\psi_j^+
-\frac{\partial\phi_j^*}{\partial\tau_a}\hat\psi_j^-]=0,
\end{equation}
$$\alpha =\overline{k+1,\infty},\quad
a=\overline{1,k},\quad j=\overline{1,\infty}$$
$$\hat H_2=
\frac{1}{2}\frac{\partial^2H}{\partial\phi_m^*\partial\phi_n^*}
\hat\psi_m^+\hat\psi_n^+
+\frac{1}{2}\frac{\partial^2H}{\partial\phi_m\partial\phi_n}
\hat\psi_m^-\hat\psi_n^-
+\frac{\partial^2H}{\partial\phi_m^*\partial\phi_n}
\hat\psi_m^+\hat\psi_n^-,\quad m,n=\overline{1,\infty}$$

Let us make use of Lemma \ref{l5}, commutation relations (\ref{e42}) and the
functional integral presentation for the scalar product \cite{B}.
It is then easy to check that if for any $z$
$$(i\frac{\partial}{\partial t}
-\frac{1}{2}\frac{\partial^2H}{\partial\phi_m^*\partial\phi_n^*}
z_mz_n
-\frac{1}{2}\frac{\partial^2H}{\partial\phi_m\partial\phi_n}
\frac{\partial}{\partial z_m}\frac{\partial}{\partial z_n}$$
$$-\frac{\partial^2H}{\partial\phi_m^*\partial\phi_n}
z_m\frac{\partial}{\partial z_n})
\int d^k\xi\exp(-\frac{1}{2}
\xi_b\xi_c
\frac{\partial\phi_j}{\partial\tau_b}
\frac{\partial\phi_j^*}{\partial\tau_c})$$
\begin{equation}\label{e43}
\exp(\xi_b
\frac{\partial\phi_j}{\partial\tau_b}z_j)
\exp(-\xi_b
\frac{\partial\phi_j^*}{\partial\tau_b}\frac{\partial}{\partial z_j})
\chi e^{\frac{1}{2}z_mM_{mn}z_n}
= 0,\qquad m,n,j=\overline{1,\infty},\qquad b,c=\overline{1,k}
\end{equation}
(where arguments $t,\tau$ at the functions $M,\chi,\phi^*,
\phi$ and arguments $\phi^*(\tau,t),\phi(\tau,t)$ at the
derivatives of $H$ are omitted) then the norm of the
vector $(i\partial/\partial t-H(\sqrt{\epsilon}\hat\psi^+,
\sqrt{\epsilon}\hat\psi^-)/\epsilon)\hat{\tilde\Psi}^{\epsilon}(t)$
tends to zero at $\epsilon\to 0$.

Formula (\ref{e43}) may be checked by calculating the Gaussian integral
at $\xi$ using presentations for $\chi, M, F, G$.

The property (\ref{e39}) is therefore proved.
Estimate now $\hat\delta^{\epsilon}(t)$. Introduce notation
$$\hat H =\frac{1}{\epsilon}H(
\sqrt{\epsilon}\hat\psi^+,
\sqrt{\epsilon}\hat\psi^-).$$
$$\hat\kappa^{\epsilon}(t)=\hat{\tilde\Psi}^{\epsilon}(t)
-\hat\Phi(t),\quad
\hat s^{\epsilon}(t) = (i\frac{\partial}{\partial t}
-\hat H)\hat{\tilde\Psi}^{\epsilon}(t).$$

The function $\hat\kappa^{\epsilon}(t)$ is the solution of the Cauchy
problem
$$i\frac{\partial}{\partial t}\hat\kappa^{\epsilon}(t)
-\hat H\hat\kappa^{\epsilon}(t) = \hat s^{\epsilon}(t),
\hat\kappa^{\epsilon}(0) = 0$$
and it is
$$\hat\kappa^{\epsilon}(t) =\int\limits_0^t
dt^{'}\exp(-i\hat H(t-t^{'}))\hat s^{\epsilon}(t^{'}).$$
Hence
$$\|\hat\kappa^{\epsilon}(t)\|\le
\int\limits_0^t\|\hat s^{\epsilon}(t^{'})\|dt^{'}
{\longrightarrow\atop \epsilon\to 0} 0$$
that implies
$$\|\hat\delta^{\epsilon}(t)\|
{\longrightarrow\atop \epsilon\to 0} 0$$

Theorem 2 is proved.

This work is supported by International Science Foundation,
grant \# MF0000.


\begin{thebibliography}{99}
\bibitem{BD}
Belov V.V., Dobrokhotov S.Yu. TMPh, 1992, v. 92, n. 2,
p. 215-254 (in Russian)

\bibitem{B}
Berezin F.A. Method of secondary quantization, M., Nauka,
1986, 318 p. (in Russian)

\bibitem{BS}
Bogoliubov N.N., Shirkov D.V. Introduction to the
Theory of Quantized Fields, - N.Y.: Interscience
Publishers, 1959

\bibitem{D}
Dirac P.A.M. Proc. Roy. Soc. A, 1927, v. 144,
pp. 234-262

\bibitem{F1}
Fock V.A. Zs. Phys., 1932, v. 75, pp. 622-647

\bibitem{F2}
Fock V.A. Soviet Phys., 1934, v. 6, p. 425

\bibitem{GJ}
Glimm J., Jaffe A. Quantum Physics.
A functional integral point of view, N.Y., 1981

\bibitem{GS}
Guillemin V., Sternberg S. Geometric asymptotics,
Providence, 1977

\bibitem{JW}
Jordan P., Wigner E.
Zh. Phys., 1928, Bd. 47, S. 631-658

\bibitem{KM}
Karasev M.V., Maslov V.P. Nonlinear Poisson brackets.
Geometry and quantization, M., Nauka, 1991, 368 p.

\bibitem{K}
Krakhnov A.D. Russ. Math. Surv., 1976, v. 31, n.3, pp. 217-218

\bibitem{LL}
Landau L.D., Livshitz E.M. Quantum mechanics. Nonrelativistic
theory, M. Nauka, 1989, 768 p. (in Russian)

\bibitem{M1}
Maslov V.P. Operational Methods, Moscow: Mir Publishers, 1976

\bibitem{M2}
Maslov V.P. Complex WKB method in nonlinear equations, M., Nauka,
1977, 384 p. (in Russian)

\bibitem{MS1}
Maslov V.P., Shvedov O.Yu. TMPh., 1994, v. 98, n. 2, pp. 266-288

\bibitem{MS2}
Maslov V.P., Shvedov O.Yu. Russian Journal of Mathematical
Physics, 1994,vol. 2, n. 2, p.217-234

\bibitem{MS3}
Maslov V.P., Shvedov O.Yu. DAN, 1994, v.338, n. 1

\bibitem{MT}
Maslov V.P., Tariverdiev S.E. Probability theory, mathematical
statistics, theoretical cybernetics, M., VINITI, 1982, v. 19,
p. 85-125

\bibitem{Sn1}
Sch\" onberg M. Nuovo cim., 1952, v. 9, n. 12, p.1139

\bibitem{Sn2}
Sch\" onberg M. Nuovo cim., 1953, v. 10, n. 4, p.419

\bibitem{Sw}
Schweber S. An introduction to relativistic quantum
field theory, Elmsford (N.Y.), 1961

\bibitem{Si}
Simon B. The $P(\Phi)_2$ Euclidean (quantum) field
theory, Princeton (N.Y.), Princeton univ. press, 1974

\end{thebibliography}
\end{document}